\documentclass[12pt,a4paper]{amsart}
\usepackage{amssymb}
\usepackage{amsmath}
\usepackage{amsthm}
\DeclareMathOperator*{\vol}{vol}

\begin{document}
\thanks{The first author was
 partially supported by NSERC grant R0010134, PSC CUNY Research Award,
 No. 60007-33-34, and Faculty Fellowship Publications Program 2003, CUNY}

\newtheorem{introtheorem}{Theorem}
\renewcommand{\theintrotheorem}{\Alph{introtheorem}}
\newtheorem{theorem}{Theorem }[section]
\newtheorem{lemma}[theorem]{Lemma}
\newtheorem{corollary}[theorem]{Corollary}
\newtheorem{proposition}[theorem]{Proposition}
\theoremstyle{definition}
\newtheorem{definition}[theorem]{Definition}
\newtheorem{example}[theorem]{Example}
\theoremstyle{remark}
\newtheorem{remark}[theorem]{Remark}

\renewcommand{\labelenumi}{(\roman{enumi})} 
\def\theenumi{\roman{enumi}}

\numberwithin{equation}{section}

\def \l {{\lambda}}
\def \a {{\alpha}}
\def \b {{\beta}}
\def \f {{\phi}}
\def \r {{\rho}}
\def \R {{\mathbb R}}
\def \H {{\mathbb H}}
\def \N {{\mathbb N}}
\def \C {{\mathbb C}}
\def \Z {{\mathbb Z}}
\def \F {{\Phi}}
\def \Q {{\mathbb Q}}
\def \e {{\epsilon }}
\def \ev {{\vec\epsilon}}
\def \ov {{\vec{0}}}
\def \GinfmodG {{\Gamma_{\!\!\infty}\!\!\setminus\!\Gamma}}
\def \GmodH {{\Gamma\backslash\H}}

\newcommand{\norm}[1]{\left\lVert #1 \right\rVert}
\newcommand{\abs}[1]{\left\lvert #1 \right\rvert}
\newcommand{\modsym}[2]{\left \langle #1,#2 \right\rangle}
\newcommand{\inprod}[2]{\left \langle #1,#2 \right\rangle}
\newcommand{\Nz}[1]{\left\lVert #1 \right\rVert_z}

\title{Modular symbols have a normal distribution}
\author{Yiannis N. Petridis}
\address{Department of Mathematics and Computer Science\\
City University of New York, Lehman College\\
250 Bedford Park Boulevard West
Bronx\\NY 10468-1589}
\email{petridis@comet.lehman.cuny.edu}
\author{Morten Skarsholm Risager}
\address{Department of Mathematical Sciences\\ University of Aarhus\\ Ny Munkegade Building 530\\ 8000 {Aa}rhus, Denmark}
\email{risager@imf.au.dk}
\date{\today}
\begin{abstract}
We prove that the modular symbols appropriately normalized and ordered have a Gaussian distribution
for all cofinite subgroups of $\hbox{SL}_2( {\mathbb R})$. We use spectral deformations to study the poles and the residues of Eisenstein series twisted by power of modular symbols.  
\end{abstract}
\maketitle

\section{Introduction}  

Let $X$ be a hyperbolic surface of finite volume. To count the number of closed prime geodesics, we 
 introduce the function $\pi (x)=\#\{c| l(c)\le x\}$, where $c$ is such a geodesic and $l(c)$ is its length.
 It follows from the
Selberg trace formula that \begin{align*}\pi (x)\sim e^x/x\end{align*} 
as $ x\rightarrow \infty$ (\cite{huber, buser}).
Error terms may be obtained but they depend on the existence of small eigenvalues of the
Laplace operator (see e.g. \cite{venkov2}). One can generalize this theorem to other
negatively curved manifolds (\cite{margulis,katok} ) and refine it as well by
counting  geodesics in a given homology class (\cite{adachisunada,phillips3}). This is achieved by
integrating the trace formula over the character variety of the surface. Since  every conjugacy class of $\pi_1(X)$
represents a free homotopy class containing a unique close geodesic, the above results can be thought of as counting 
group elements in $\pi_1(X)$ with weight $1$ and with weight  the characterstic
function of the homology class, respectively.

In this paper we use weights which are polynomials in the values of the Poincar\'e pairing.
We have the  pairing between homology and
cohomology: 
\begin{align*}H_1(X,\R)\times H^1_{\rm dR}(X,\mathbb R)\rightarrow \mathbb R\end{align*} and a map
 $\phi:\Gamma =\pi_1(X)\rightarrow H_1(X,\Z)$. Let $\modsym{\cdot}{\cdot}$ be the
composition of the two maps. Using moments we examine the distribution of the values of the composition for a fixed cohomology class. We restrict ourselves to surfaces with cusps 
and cuspidal cohomology. We fix a cuspidal cohomology class $\a$, which
we can take to be harmonic. We have no loss of generality by assuming that $\alpha =\Re (f(z)dz)$, where $f(z)$ is a holomorphic cusp form of weight $2$. 
The group $\Gamma$ can be realized as a subgroup of $\hbox{SL}_2( {\mathbb R})$. 
Let an element $\gamma \in \Gamma$ have lower row $(c, d)$.
We use the normalization
$$\modsym{\gamma}{\a}=-2\pi i \int_{\phi (\gamma)}\alpha.$$
\begin{introtheorem}\label{maintheoremforreal}    The values $\modsym{\gamma}{\a}$ appropriately
normalized have a
normal distribution. More precisely 
\begin{equation}\frac{ \#\left\{\gamma\in (\GinfmodG)^T \left| \frac{\widetilde{\langle \gamma ,\alpha\rangle }}{i\sqrt{\log (c^2+d^2)}}\in [a,b]\right. \right\} }{\#(\GinfmodG)^T}\to \frac{1}{\sqrt{2\pi}}\int_a^b\!\! \exp\left(-\frac{x^2}{2}\right)dx\end{equation} as $T\to\infty$.
\end{introtheorem}

Here   $(\GinfmodG)^T$ is set of  elements in  $\GinfmodG$ with $c^2+d^2\leq T$ while
\begin{align*}\widetilde{\modsym{\gamma}{\a}}=\sqrt{\frac{\vol{(\GmodH)}}{{8\pi^2\norm{f}^2}}}\modsym{\gamma}{\a},\end{align*}
where $\norm{f}$ is the Petersson norm of $f$.
In fact we can consider complex valued 1-forms $f(z)dz$, where $f$ is a
holomorphic cusp form of weight 2. Then

\begin{introtheorem}\label{maintheorem}
Asymptotically  $\frac{\widetilde{<\gamma,f>}}{\sqrt{\log (c^2+d^2)}}$ has bivariate Gaussian distribution with correlation coefficient zero. More precisely we have for $R\subset\C$ a rectangle
\begin{equation*}\frac{ \#\left\{\gamma\in (\GinfmodG)^T \left| \frac{\widetilde{\langle \gamma,f\rangle }}{\sqrt{\log (c^2+d^2)}}\in R\right. \right\} }{\#(\GinfmodG)^T}\to \frac{1}{{2\pi}}\int_R\!\! \exp\left(-\frac{x^2+y^2}{2}\right)dxdy\end{equation*} as $T\to\infty$.
\end{introtheorem}
Here \begin{align*}\widetilde{\modsym{\gamma}{f}}=\sqrt{\frac{\vol{(\GmodH)}}{{8\pi^2\norm{f}^2}}}\modsym{\gamma}{f}.\end{align*}
This work uses heavily Eisenstein series twisted by modular symbols,
introduced by Goldfeld. The general framework is as follows. Let $f(z), g(z)$ be  holomorphic cusp forms of weight $2$ for a fixed cofinite discrete
subgroup $\Gamma$ of $\hbox{SL}_2({\mathbb R})$. 
In \cite{goldfeld1, goldfeld2} Goldfeld introduced Eisenstein series 
associated with modular symbols  defined in a right half-plane as
\begin{equation}\label{twistedseries}E^{m, n}(z, s)=\sum_{\gamma \in\GinfmodG}
\modsym{\gamma}{f}^m\overline{\modsym{\gamma}{g}}^n \Im (\gamma z)^s,
\end{equation}
where for $\gamma\in \Gamma$ the modular symbol $\modsym{\gamma}{f}$ is given by
\begin{equation}
\modsym{\gamma}{f} =-2\pi i \int_{z_0}^{\gamma z_0}f(z)\, dz ,
\end{equation}
and one defines similarly $\modsym{\gamma}{g}$.
Here $z_0$ is an arbitrary point in the upper half-plane $\H$.

If we take $f(z)$ to be a Hecke eigenform for $\Gamma_0(N)$ with rational coefficients and $E_f$ is 
the elliptic curve over $\mathbf Q$ corresponding to it by the Eichler-Shimura
 theory, then
$$\modsym{\gamma}{f}=n_1(f, \gamma )\Omega_1(f)+n_2(f, \gamma )
\Omega_2(f),$$
where $n_i\in\mathbb Z$ and $\Omega_i$ are the periods of $E_f$.
The conjecture $n_i\ll N^k$ for $\abs{c}\le N^2$ and some fixed  $k$ (Goldfeld's 
conjecture) is equivalent
to Szpiro's conjecture $D\ll N^C$ for some $C$, where $D$ is the discriminant
 of $E_f$. This was the  motivation to study the distribution of 
modular symbols.

As an example of such a distributional result Goldfeld conjectured in \cite{goldfeld1} that 
\begin{equation}\label{hochstapler}
\sum_{c^2+d^2\le T}\langle \gamma ,f\rangle\sim R(i)T,
\end{equation}
where $R(z)$ is the residue at $s=1$ of $E^{1, 0}(z, s)$, and we sum over the
 elements in $\GinfmodG$ with lower row $(c, d)$. 
This in now proved in \cite[Theorem 7.3]{gos}.
He also suggested that, when $f=g$, the twisted Eisenstein series $E^{1, 1}(z, s)$ has a simple pole at $s=1$ with the zero Fourier coefficient of the 
residue proportional to the Petersson
norm $\norm{f}^2$. He concludes the conjectural asymptotic formula
\begin{equation}\label{hochstapler2}
\sum_{c^2+d^2\le T}\abs{\langle \gamma, f\rangle}^2\sim R^*(i)T,
\end{equation} where $R^*(z)$ is the residue of $E^{1, 1}(z, s)$ at $s=1$ and where the summation is again over matrices in $\GinfmodG$  with lower row $(c, d)$.
 In this work we, among other things, reprove (\ref{hochstapler}) while settling (\ref{hochstapler2}) in the negative. But our result shows that the Petersson norm does indeed play a role, see Theorem \ref{eswareinkoeniginthule} below.
Averages of functions of modular symbols have been investigated also in \cite{manin}.

It turns out to be crucial to consider Eisenstein series
 associated with the real harmonic differentials $\a_i=\Re(f_i(z)dz)$
 or $\a_i=\Im(f_i(z)dz)$   where $f_i$ are  holomorphic cusp forms of 
weight two. We shall write \begin{equation} \modsym{\gamma}{\a_i}=-2\pi i\int_{z_0}^{\gamma z_0}\a_i .\end{equation}
As in \cite{petridis}
we define 
\begin{equation}\label{stpraxedis}
E(z, s,\ev)=\sum_{\gamma\in\GinfmodG}                    
\chi_{\ev}(\gamma ) \Im (\gamma z)^s,
\end{equation}
where $\chi_{\ev}$ is an $n$-parameter family of characters of the group
defined by 
\begin{equation}
\chi_{\ev}(\gamma )=\exp \left(-2\pi i\left( \sum_{k=1}^n\e_k\int_{z_0}^{\gamma z_0}\alpha_k \right) \right).
\end{equation}
The convergence of this is guaranteed for $\Re (s)>1$ by comparison with the
 standard Eisenstein series. 
The Eisenstein series with a  character transform as 
\begin{equation}
E(\gamma z, s,\ev)={\bar\chi_{\ev }(\gamma )}E(z, s,\ev).
\end{equation}
In the domain of absolute convergence  we see that 
\begin{equation}\label{vorueber}
\left.\frac{\partial^{n}E(z, s,\ev)}{\partial\epsilon_1\ldots\partial \epsilon_n}\right|_{\ev=\ov}=\sum_{\gamma\in\GinfmodG}\prod_{i=1}^n\modsym{\gamma}{\alpha_i}\Im(\gamma z)^s,
\end{equation}
by termwise differentiation. By taking linear combinations of these we may of course recover the original series (\ref{twistedseries}).
This observation allowed the first author to give a new approach to the Eisenstein
series twisted with modular symbols using perturbation theory. In particular, 
a new proof of the analytic continuation was given in \cite{petridis}
and the residues of $E^{1, 0}(z, s)$ on the critical line were identified.
In this paper we further pursue this method. We start by giving a third much shorter proof of the main theorem in \cite{osullivan}.
\begin{introtheorem}[\cite{osullivan,petridis}]\label{ganztreubisandasgrab}
The functions $E^{m,n}(z,s)$ have meromorphic continuation to the whole $s$-plane. In $\Re (s)>1$ the series are absolutely  convergent and, consequently, they are analytic.
\end{introtheorem} 
The last claim of the theorem is new and enables us to evaluate the growth of the modular symbols as $\gamma$ runs through the group $\Gamma$. The best known result in this aspect is 
 \begin{equation*}\modsym{\gamma}{f}=O(\log(c^2+d^2)).\end{equation*}
 This is due to Eichler (see \cite{eichler}). Using the above theorem we get the following slightly weaker result.
\begin{introtheorem} \label{sonnenallee}For any $\varepsilon>0$ we have
 \begin{equation*}\modsym{\gamma}{f}=O_{\varepsilon}((c^2+d^2)^{\varepsilon}).\end{equation*}
\end{introtheorem}

We then continue to study the singularity of $E^{m,n}(z,s)$ at $s=1$ when $f=g$. In particular we study the pole order
 and the leading term in the singular part of the Laurent expansion. In principle the method gives 
the full Laurent expansion of $E^{m,n}(z,s)$ but only in terms of the coefficients in the
  Laurent expansions of the resolvent kernel and the usual nonholomorphic Eisenstein series at $s=1$.
 The combinatorics involved in getting useful expressions are quite ponderous. As a result we settle with calculating some of the most interesting  coefficients and evaluate the pole orders.  As an example of this type of result we  have:
\begin{introtheorem}\label{yetanotherstupidlabel} At $s=1$, $E^{2,0}(z,s)$ has a simple pole with residue $$\frac{1}
{{\vol({\GmodH})}}\left(2\pi i\int_{i\infty}^zf(z)dz\right)^2$$ 
while $E^{1,1}(z, s)$ has a double pole with residue
$$\frac{4\pi^2}{{\vol({\GmodH})}}\abs{\int_{i\infty}^zf(z)dz}^2+\frac{16\pi^2}{\vol ({\GmodH})}
\int_{\GmodH}(E_0(z')-r_0(z,z')){y'}^2\abs{f(z')}^2d\mu(z').$$ 
The coefficient of $(s-1)^{-2}$ is
$$ \frac{16 \pi^2 \norm{f}^2}{\vol(\GmodH)^2}.$$
\end{introtheorem}
Here  the coefficient $r_0(z,z')$  is the constant term in the Laurent expansion of the resolvent kernel around $s=1$.
The coefficient $E_0(z)$ is the constant term in the Laurent expansion of the usual nonholomorphic Eisenstein series and is given by Kronecker's limit formula. 
For $\Gamma=\hbox{SL}_2(\Z )$ this is classical, see, for instance 
\cite[p.~273--275]{lang}. For a generalization to all $\Gamma$ see \cite{goldstein}.

We wish to use these results to obtain results \`a la (\ref{hochstapler}). We do this using the method of contour integration but, in order to make this work, we need to prove a result on the growth of $E^{m,n}(z,s)$ as $\Im(s)\to\infty$. We can prove
\begin{introtheorem}\label{demsterbendseinebuhle} The functions $E^{m,n}(z,s)$ grow at most polynomially on vertical lines with $\sigma >1/2$. More precisely: for every $\varepsilon >0$ and $\sigma \in (1/2, 1]$ and $z\in K$, a compact set, we have
\begin{equation}
E^{m, n}(z, \sigma+it)=O(|t|^{(6(m+n)-1)(1-\sigma)+\varepsilon}).
\end{equation}
\end{introtheorem} 
Using the above theorems and contour integration  we get asymptotic expansions for summatory functions like  the one in  (\ref{hochstapler}). An example of the results we prove is:
\begin{introtheorem}\label{eswareinkoeniginthule}
There exists $\delta>0$ such that
\begin{align*}
\sum_{c^2+d^2\leq T}\modsym{\gamma}{f}^2&=\frac{1}{\vol(\GmodH)}\left(-2\pi i\int_{i\infty}^zf(\tau )\, d\tau\right)^2T+O(T^{1-\delta})\\
\sum_{c^2+d^2\leq T}\abs{\modsym{\gamma}{f}}^{2}&=\frac{(16\pi^2)}{\vol(\GmodH)^2}\norm{f}^{2}T\log T+O(T). 
\end{align*}
The summations are over $(c, d)$ lower row of $\gamma\in \GinfmodG$.
\end{introtheorem}
This settles the status of (\ref{hochstapler2}). How small we can make $1-\delta$ in the above theorem depends on how good polynomial
 bounds we have in Theorem \ref{demsterbendseinebuhle} and whether the Laplacian has small eigenvalues. If there are no such eigenvalues we can prove 
\begin{align*}1-\delta=\frac{12}{13}+\varepsilon.\end{align*}
By using similar asymptotic expansions  we can calculate the  moments of the normalized modular symbols 
and  prove the distributional result in Theorem \ref{maintheorem},  which is the main theorem of our work.

The idea of putting the Eisenstein series in a continuous family to study how the spectrum changes as the parameters change is very fruitful, see for instance \cite{bruggeman}. In fact it is possible to construct a proof of the first part of Theorem \ref{ganztreubisandasgrab} different from the proof given in this paper using ideas in \cite[Chapter 15]{bruggeman}. 

The study of $E^{m, n}(z, s)$ using perturbed Eisenstein series is an
interesting  application
of the spectral deformations used in \cite{phillipssarnak1, phillipssarnak2, phillips3, petr}. Our contribution in \cite{petridis} was to put the
 Eisenstein series
with modular symbols into this framework.
In this work we apply the same techniques to produce results which at least to us seem difficult to attack with the methods used by Goldfeld,  O'Sullivan et.al.

\section{Finding Laurent expansions using perturbation theory}\label{polarparts}

We assume that $\Gamma$ has only one cusp and that this is of width 1. The generalization to the multiple cusp case is straightforward. We note that we can always assume that $\alpha_i$ is the real part of a holomorphic cusp form since 
\begin{equation*}\Im(f(z)dz)=\Re(-if(z)dz)\end{equation*} and $-if$ is a holomorphic cusp form of weight two.
 We want to approximate the real differentials $\alpha_i$ with compactly supported ones. We do this as follows. 
Let  
\begin{equation*}f(z)=\sum_{n=1}^\infty a_ne^{2\pi i n z}\end{equation*} 
be the Fourier expansion of $f(z)$. We  define 
\begin{equation*}F(z)=\sum_{n=1}^\infty \frac{a_n}{2\pi i n}e^{2\pi i n z}.\end{equation*} If $\alpha_i(z)=\Re(f(z)dz)$ then also $\alpha_i(z)=d\Re(F(z))$. We choose a fundamental domain $F$ in such a way that there exists $T_0\in\R_+$, $a,b\in\R$ such that 
\begin{equation*}F\cap\{z\in\H\, \,|\Im(z)>T_0\}=\{z\in\H\,\,|\Im(z)>T_0,\quad 0< \Re(z)\leq 1\}.
\end{equation*} We then choose a smooth function $\tilde\psi:\R\to[0,1]$ such that 
\begin{equation*}\tilde \psi(t)=\begin{cases}1,&\textrm{if }t\leq 0, \\0, &\textrm{if }t\geq 1.\end{cases}\end{equation*}We then define, for $T>T_0$, $\psi^T:F\to[0,1]$ by $\psi^T(z)=\tilde\psi(\Im(z)-T)$. For $z\in F$ we define \begin{eqnarray*}w_i^T&=&d(\psi^T\Re(F))\\g_i^T&=&(1-\psi^T)\Re(F),\end{eqnarray*} and extend these to smooth $\Gamma$ automorphic 1-forms on  $\H$ by setting $w_i^T(\gamma z)=w_i^T(z)$ and $g_i^T(\gamma z)=g_i^T(z)$ for each $\gamma\in\Gamma$. Then we have \begin{equation}\alpha_i=w_i^T+dg_i^T.\end{equation}
The following proposition is easy to verify. 
\begin{proposition}\label{approxprop}\quad
\begin{enumerate}
\item The smooth 1-form $w_i^T$ is compactly supported on $\GmodH$.
\item \label{espresso} If $z\in F$ and $\Im(z)\leq T$ then $w_i^T(z)=\alpha_i(z)$.
\item \label{valrohna}If $z\in F$ and $\Im(z)\leq T$ then $\int_{i\infty}^z w_i^T=\int_{i\infty}^z\alpha_i$.
\item \label{thisisfun} $\modsym{\gamma}{\alpha_i}=\modsym{\gamma}{w_i^T}$ for all $\gamma\in\Gamma$ and all $T> T_0$.
\end{enumerate}
\end{proposition}
We shall often exclude $T$ from the notation and simply write \begin{equation*}\alpha_i=w_i+dg_i.
\end{equation*}
We note that by Proposition \ref{approxprop} (\ref{thisisfun})
\begin{equation}
\chi_{\ev}(\gamma )=\exp \left(-2\pi i\left(
\sum_{k=1}^{n}\e_k\int_{z_0}^{\gamma z_0} w_k \right) \right).
\end{equation}
 We consider the space $L^2(\GmodH , \bar\chi_{\ev})$ of square integrable functions which transform as 
$$h(\gamma \cdot z)=\bar\chi_{\ev}(\gamma )h(z), \quad \gamma\in\Gamma$$
under the action of the group. We introduce  unitary operators
$$U(\ev):L^2(\GmodH )\to L^2(\GmodH ,
\bar\chi_{\ev})$$
given by
$$(U(\ev)h)  (z):=U(z,\ev)h(z)=\exp\left(2\pi i \left(\sum_{k=1}^{n}\e_k \int_{i\infty}^zw_k \right)\right)h(z).$$
We set 
$$L(\ev)=U(\ev)^{-1}\Delta U(\ev)$$
and
\begin{equation}\label{girlwithapearlearring}
E(z, s,\ev)=U(\ev)D(z, s,\ev).
\end{equation}
We note that $L(\ev)=\Delta$ {\lq close to the cusp\rq} since $U(\ev)$ is compactly supported.
We note also that $E(z,s,\ev)$ is independent of the choice of differential within a cohomology class, i.e. independent of $T$,while $D(z,s,\ev)$ and $U(\ev)$ are not.  We also remark that \cite[Remark 2.2]{petridis} is only true for $z_0=i\infty$, since both $E(z, s,\ev)$ and $D(z, s,\ev)$ have asymptotic behavior
at $\infty$ of the form $y^s$ for $\Re (s)>1$ and, consequently,
$U(z,\epsilon)$ should tend to $1$, as $\Im(z)\to\infty$.
We define  $\langle f_1 dz+f_2d\bar z, g_1dz+g_2d\bar z\rangle =2y^2(f_1\bar g_1+f_2\bar g_2)$, 
$\delta (pdx+qdy)=
-y^2(p_x+q_y)$.
\begin{lemma}\label{hetstraatje}
The conjugated operator $L(\ev) $ is given by
\begin{eqnarray}\label{viewofdelft}
L(\ev)h&=&\Delta h+4\pi i \sum_{k=1}^{n}\epsilon_k\langle dh, w_k\rangle-2\pi i \left(\sum_{k=1}^{n}\e_k \delta(w_k)\right) h
\nonumber \\ &&-4\pi ^2 \left(\sum_{k,l=1}^{n}\e_k\e_l\langle w_k,w_l\rangle\right)h.
\end{eqnarray}
\end{lemma}
\begin{proof}
The proof uses induction on $n$. The result for $n=1$ may be found in \cite[p.~113]{petr}. With the convention that $U(\e_k)=U((0,\ldots,0,\e_k, 0,\ldots,0))$ we see that
\begin{eqnarray*}L(\ev)h&=&U(\e_n)^{-1}U(\e_1,\ldots\e_{n-1},0)^{-1}\Delta U(\e_1,\ldots\e_{n-1},0) U(\e_n)h\\&=&U(\e_n)^{-1}\left(\Delta U(\e_n)h+4\pi i \sum_{k=1}^{n-1}\epsilon_k\langle dU(\e_n)h, w_k\rangle\right.\\&&\left.-2\pi i \left(\sum_{k=1}^{n-1}\e_k \delta(w_k)\right) U(\e_n)h
    -4\pi ^2 \left(\sum_{k,l=1}^{n-1}\e_k\e_l\langle w_k,w_l\rangle\right) U(\e_n)h \right).
\end{eqnarray*}
We apply the result for one variable once more in the $\e_n$ variable and use the
chain rule in the form
$$d(U(\e_n)h)=U(\e_n)dh+2\pi i \e_nU(\e_n)hw_n $$
to get the result.
\end{proof}
 In the rest of the paper we will use the following convention. A function with a subscript variable will denote the partial derivative of the function in that variable.
Lemma \ref{hetstraatje} gives
\begin{equation}\label{Lderivedonce}
L_{\e_k}(\ov)h=4\pi i \langle dh, w_k\rangle-2\pi i (\delta w_k) h,
\end{equation}
\begin {equation}\label{Lderivedtwice}
L_{\e_k\e_l}(\ov)h=-8\pi ^2\langle w_k, w_l\rangle h .
\end{equation}
and all higher order derivatives vanish. Differentiating the eigenvalue equation (see \cite[Lemma 8]{petridis})  
\begin{equation}\label{eigenvalueequation}(L(\ev)+s(1-s))D(z, s,\ev)=0\end{equation}
and applying the resolvent of the Laplace operator, $R(s)=(\Delta +s(1-s))^{-1}$,  we get
\begin{equation}\label{denmarkiscold}
D_{\epsilon_k}(z, s,\ov)=-R(s)\left(L_{\e_k}(\ov)D(z,s,\ov)\right)
\end{equation}
and
\begin{eqnarray}\label{denmarkisrainfull}
D_{\epsilon_1,\ldots,\epsilon_n}(z, s,\ov)&=-R(s)&\left(\sum_{k=1}^{n}L_{\e_k}(\ov)D_{\e_1,.,\widehat{\e_k},.,\e_n}(z,s,\ov)\right.\\\nonumber&&\left.+\sum_{\substack{k,l=1\\ k<l}}^{n}L_{\e_k\e_l}(\ov)D_{\e_1,.,\widehat{\e_k},.,\widehat{\e_l},.,\e_n}(z,s,\ov)\right).
\end{eqnarray}
 Here $\widehat{\e_k}$ means that we have excluded $\e_k$ from the list. The validity of the inversion of $(\Delta+s(1-s))$ using the resolvent follows from the following lemma:
\begin{lemma}\label{squareintegrable} Let $n\geq 1$. For $\Re(s)$ sufficiently large we have 
  \begin{equation*}
    D_{\e_1,\ldots,\e_n}(z,s,\ov)\in L^2(\GmodH,d\mu).
  \end{equation*}
\end{lemma}
\begin{proof}
Since the function $D(z,s,\ev)$ is $\Gamma$-automorphic we see that also  $D_{\e_1,\ldots,\e_n}(z,s,\ov)$ is $\Gamma$-automorphic. From (\ref{girlwithapearlearring}) we obtain  that 
\begin{equation}D_{\e_1,\ldots,\e_n}(z,s,\ov)=\!\!\!\!\sum_{\!\!\!\vec m\in\{0,1\}^n}\prod_{k=1}^n\left(-2\pi i\int_{i\infty}^zw_k\right)^{m_k}E_{\e_1^{1-m_1},\ldots,\e_n^{1-m_n}}(z,s,\ov).\end{equation} 
We note that since $w_i$ is compactly supported all the terms with $\vec m\neq \ov$ becomes compactly supported. Now in order to control the term with $\vec m=\ov$ we need some bound on the growth of $\modsym{\gamma}{\a_i}$. Any bound of the form 
$$\abs{\modsym{\gamma}{\alpha_i}}\leq C (c^2+d^2)^b$$ will do. We quote \cite[Lemma 1.1]{osullivan} with $z=i$ to get $b=1$. If we use the inequality 
\begin{equation}\label{yetanotherstupidbound}(c^2+d^2)\leq \frac{|cz+d|^2}{y}\frac{1+|z|^2}{y}\end{equation}
which follows from adding $\abs{cz+d}^2\geq (cy)^2$ and $\abs{z}^2\abs{cz+d}^2\geq (dy)^2$ (see also \cite[Lemma 4]{knopp}) we get:
\begin{align*}
  \abs{E_{\e_1,\cdots,\e_n}(z,s,\ov)}&\leq \sum_{\substack{\gamma \in\GinfmodG\\\gamma\neq I}}\abs{\prod_{j=1}^n\modsym{\gamma}{\a_j}}\Im(\gamma z)^\sigma\\
&=C\left(\frac{1+\abs{z}^2}{y}\right)^n \sum_{\substack{\gamma\in\GinfmodG\\\gamma\neq I}}\Im(\gamma z)^{\sigma-n}
\intertext{We note that the sum is $O_\sigma(y^{1-\sigma+n})$ by \cite[p. 13]{kubota} so we get}
&\leq C'\left(\frac{1+\abs{z}^2}{y}\right)^n y^{1-\sigma+n}\\
&\leq C'' y^{1-\sigma+2n}
\end{align*}
Hence we conclude that for $\sigma>2+2n$ we have  $D_{\e_1,\ldots,\e_n}(z,s,\ov)\in L^2(\GmodH,d\mu(z))$.  
\end{proof}

 Using the representation (\ref{denmarkisrainfull}) we may give a short proof of the analytic continuation of the functions defined in a half-plane by (\ref{twistedseries})
\begin{lemma}The functions $D_{\e_1\ldots\e_n}(z, s,\ov)$ have meromorphic continuation to $\C$. In $\Re(s)>1$ the functions are analytic.
\end{lemma}
\begin{proof}The proof uses induction on $n$. For $n=0$ the function is the classical Eisenstein series and one of the many known proofs may be found in \cite{kubota}. We note that  by (\ref{Lderivedonce}) and (\ref{Lderivedtwice}) $L_{\e_k}(\ov)D_{\e_1,.,\widehat{\e_k},.,\e_n}(z,s,\ov)$ and $L_{\e_k\e_l}(\ov)D_{\e_1,.,\widehat{\e_k},.,\widehat{\e_l},.,\e_n}(z,s,\ov)$ are compactly supported. Hence from (\ref{denmarkisrainfull}) and \cite[Theorem 1]{mueller2} the conclusion follows.  
\end{proof}
\begin{remark}\label{karlwantsshoes} From the above lemma, (\ref{girlwithapearlearring}) and (\ref{vorueber}) we find that  $$\left.\frac{\partial^{n}E(z, s,\ev)}{\partial\epsilon_1\ldots\partial \epsilon_n}\right|_{\ev=\ov}$$ has meromorphic continuation and that in  $\Re (s)>1$ these functions are analytic. By taking linear combinations of these (see (\ref{vorueber})) we obtain the Theorem \ref{ganztreubisandasgrab}.\end{remark}
\begin{proposition}\label{absconvergence}
  The sum defining $E^{m,n}(z,s)$ is absolutely convergent whenever $\Re(s)>1$.
\end{proposition}
\begin{proof}
  Note that if we can prove the above for $f=g$ and $m=n$ then we get the general result by appealing to the elementary inequality $2ab\leq a^2+b^2$ for $a,b\in\R$. This gives $$\abs{\langle \gamma,f \rangle ^m\overline{\langle\gamma,g\rangle^n}}\leq\frac{1}{2}\left(\abs{\langle\gamma,f\rangle}^{2m}+\abs{\langle\gamma,g\rangle}^{2n}\right),$$ and comparison with  $f=g$ and $m=n$ type Eisenstein series gives the result.

For the case  $f=g$ and $m=n$, the proof uses Landau's result. He proved that Dirichlet series with positive coefficients has a singularity on the line of absolute convergence, see e.g. \cite[Section 9.2]{titchmarsh}.
By Remark \ref{karlwantsshoes} we get the first singularity of $E^{m,m}(z, s)$ at $s=1$ or further to the left.
\end{proof}
Clearly $E_{\e_1,\ldots,\e_n}(z,s,\ov)$ is also absolutely convergent for $\Re(s)>1$ by the same proof.  
We immediately get the following corollary:
\begin{corollary}\label{slowgrowth}For any fixed $z\in \H$, $\varepsilon>0$ we have 
\begin{align*}
\modsym{\gamma}{f}&=o(|cz+d|^\varepsilon)\\ \modsym{\gamma}{\alpha}&=o(|cz+d|^\varepsilon) 
\end{align*} as $|cz+d|\to \infty$. 
\end{corollary} 
\begin{proof} Since the terms in an absolutely convergent series tend to zero Proposition \ref{absconvergence} implies that for any $m\in N$,
$$\modsym{\gamma}{f}^m\Im(\gamma z)^{2}=\modsym{\gamma}{f}^m\frac{y^2}{\abs{cz+d}^{4}}\to 0.$$ Hence $\modsym{\gamma}{f}=o(|cz+d|^{4/m})$. Similarly with $\modsym{\gamma}{\alpha}$.
\end{proof}
We note that by picking $z=i$ we get Theorem \ref{sonnenallee}.
\begin{remark}We note that since 
$$D(z,s,\ev)=U(-\ev)E(z,s,\ev)=\sum_{\gamma\in\GinfmodG}\exp\left(\sum_{k=1}^n-2\pi i\e_k \int_{i\infty}^{\gamma z}w_k\right)\Im(\gamma z)^s, $$  we find that 
the function $D_{\e_1,\ldots\e_n}(z,s,\ov)$ has the series representation
\begin{equation}\label{seriesrep}D_{\e_1,\ldots\e_n}(z,s,\ov)=\sum_{\gamma\in\GinfmodG}\prod_{k=1}^{n}\left(-2\pi i \int_{i\infty}^{\gamma z}w_k\right)\Im(\gamma z)^s,\end{equation} whenever the series is convergent. We also remark that from (\ref{girlwithapearlearring}) we have
\begin{equation}\label{DtoE}E_{\e_1,\ldots,\e_n}(z,s,\ov)=\!\!\!\!\sum_{\!\!\!\vec m\in\{0,1\}^n}\prod_{k=1}^n\left(2\pi i\int_{i\infty}^zw_k\right)^{m_k}D_{\e_1^{1-m_1},\ldots,\e_n^{1-m_n}}(z,s,\ov).\end{equation} 
\end{remark}
 Combining this with Proposition \ref{approxprop} (\ref{valrohna}) and (\ref{thisisfun}),  we see that, if $z\in F$ and $\Im(z)<T$, then 
\begin{equation}\label{getit}D_{\e_1,\ldots\e_n}(z,s,\ov)=\sum_{\gamma\in\GinfmodG}\prod_{k=1}^{n}\left(-2\pi i \int_{i\infty}^{\gamma z}\a_k\right)\Im(\gamma z)^s.\end{equation} 
In particular 
\begin{equation}\label{coffeeisgood}\lim_{T\to\infty}D_{\e_1,\ldots\e_n}(z,s,\ov)=\sum_{\gamma\in\GinfmodG}\prod_{\substack{k=1}}^{n}\left(-2\pi i \int_{i\infty}^{\gamma z}\a_k\right)\Im(\gamma z)^s,\end{equation} for all $z\in\H$.

\begin{lemma}\label{stupidbound} For $\sigma>1$ we have  
\begin{equation}\sum_{\gamma\in\GinfmodG}\abs{\prod_{k=1}^{n}\left(-2\pi i \int_{i\infty}^{\gamma z}\a_k\right)}\Im(\gamma z)^\sigma=O(y^{1-\sigma+\e})\end{equation}
as $\Im(z)\to \infty$ for  $z\in F$. In particular $\lim_{T\to\infty}D_{\e_1,\ldots\e_n}(z,\sigma+it,\ov)=O(y^{1-\sigma})$.\end{lemma}
\begin{proof}
We have for $\sigma>1$ (see \cite[p.13]{kubota})\begin{equation}\label{Ebound}\sum_{\substack{\gamma\in\GinfmodG\\\gamma\neq I}}\Im(\gamma z)^\sigma=O_\sigma(y^{1-\sigma})\end{equation} as $\Im(z)\to\infty$. From Corollary \ref{slowgrowth} we see that if we fix e.g. $z_0=i$ there exists a constant $C>0$ such that
\begin{equation*}\abs{\prod_{k=1}^{n}\modsym{\gamma}{\a_k}}\leq C \Im(\gamma z_0)^{-\varepsilon}.\end{equation*}  this gives, using $\modsym{I}{\a_k}=0$,
\begin{equation*}\sum_{\gamma\in\GinfmodG}\abs{\prod_{k=1}^{n}
\modsym{\gamma}{\a_k}}\Im(\gamma z)^\sigma\leq C\sum_{\substack{\gamma\in\GinfmodG\\\gamma\neq I}}\Im(\gamma i)^{-\varepsilon}\Im(\gamma z)^\sigma\end{equation*} If we use the inequality (\ref{yetanotherstupidbound}) this is majorized by 
\begin{equation*}
  C\sum_{\substack{\gamma\in\GinfmodG\\\gamma\neq I}}\Im(\gamma z)^{\sigma-\varepsilon}\left(\frac{1+|z|^2}{y}\right)^\varepsilon=O_\sigma(y^{1-\sigma+\e}).
\end{equation*} In the last equality we used (\ref{Ebound}).
The claim now follows by induction from (\ref{DtoE}) by isolating $D_{\e_1,\ldots,\e_n}(z,s,\ov)$, using (\ref{seriesrep}) and the fact that 
\begin{equation*}-2\pi i \int_{i\infty}^z\a_k\end{equation*} is $O(e^{-2\pi y})$ as $\Im(z)\to\infty.$
\end{proof}

\begin{lemma}\label{anotherstupidbound}
For $\Re(s)>1$ we have
\begin{equation}\nonumber\int_{\GmodH}|\langle d\lim_{T\to\infty} D_{\e_1,\ldots,\hat\e_j,\ldots\e_n}(z,s,\ov),\alpha_{j}\rangle| d\mu(z)<\infty .
\end{equation}
\end{lemma}
\begin{proof} 

 Using (\ref{coffeeisgood}) we see that for $\Re(s)>1$ 
\begin{equation}
d\lim_{T\to\infty}D_{\e_1,\ldots,\hat\e_j,\ldots,\e_n}(z,s,\ov)=\sum_{\gamma\in\GinfmodG}d\left(\prod_{\substack{k=1\\k\neq j}}^{n}\left(-2\pi i \int_{i\infty}^{\gamma z}\a_k\right)\Im(\gamma z)^s\right).
\end{equation}
Using\begin{eqnarray*}d\left(-2\pi i\int_{i\infty}^{\gamma z}\a_k\right)&=&-2\pi i\a_k\\ 
d\Im(\gamma z)^s&=&\frac{s}{2y}\left(-i\left(\frac{c\overline z+d}{cz+d}\right)\Im(\gamma z)^sdz+i\left(\frac{c z+d}{c\overline z+d}\right)\Im(\gamma z)^sd\overline z\right)\end{eqnarray*} we find that 
\begin{align*}
\langle d\lim_{T\to\infty}D_{\e_1,\ldots,\hat\e_j,\ldots,\e_n}(z,s,\ov)&,\a_j \rangle=\\2y^2\left[\sum_{\substack{l=1\\l\neq j}}^n\right.&\left(-2\pi i \frac{f_l}{2}\frac{\overline{f_j}}{2}\sum_{\gamma\in\GinfmodG} \prod_{\substack{k=1\\k\neq j,l}}^{n}\left(-2\pi i \int_{i\infty}^{\gamma z}\a_k\right)\Im(\gamma z)^s\right)\\
+\frac{-is}{2y}&\left.\frac{\overline{f_j}}{2}\sum_{\gamma\in\GinfmodG} \prod_{\substack{k=1\\k\neq j}}^{n}\left(-2\pi i \int_{i\infty}^{\gamma z}\a_k\right)\left(\frac{c\overline z+d}{cz+d}\right)\Im(\gamma z)^s\right]\\
&+ \textrm{complex conjugate.}
\end{align*}
The claim now follows from Lemma \ref{stupidbound}, since $f_i(z)=O(e^{-2\pi y})$ as $\Im(z)\to \infty$.
\end{proof}
Using this lemma we can prove the following important result
\begin{lemma}\label{almostregular}
For all $j=1,\ldots,n$ and $\Re (s)>1$
  \begin{equation*}\int_{\GmodH}\langle dD_{\e_1,\ldots,{\widehat{ \e_j}} ,\ldots\e_n}(z,s,\ov),w_{j}\rangle d\mu(z)\to 0 \textrm{ as }T\to\infty .\end{equation*}
 \end{lemma}
\begin{proof}We start by showing that
 \begin{align}
\nonumber\int_{\GmodH}&\langle d\lim_{T\to\infty} D_{\e_1,\ldots,\hat\e_j,\ldots\e_n}(z,s,\ov),\alpha_{j}\rangle d\mu(z)=0.
\end{align}
If we let $F_M=\{z\in F| \Im(z)\leq M\}$ then by lemma \ref{anotherstupidbound} the left-hand side is
 \begin{align*}\nonumber\int_{F_M}&\langle d\lim_{T\to\infty} D_{\e_1,\ldots,\hat\e_j,\ldots\e_n}(z,s,\ov),\alpha_{j}\rangle d\mu(z)+\e(M) \end{align*}
where $\e(M)\to 0$ as $M\to \infty$. We have
\begin{align}\nonumber\int_{F_M}&\langle d\lim_{T\to\infty} D_{\e_1,\ldots,\hat\e_j,\ldots\e_n}(z,s,\ov),\alpha_{j}\rangle d\mu(z)
=\\&\int_{F_M}\frac{\partial}{\partial z}\left(\lim_{T\to\infty} D_{\e_1,\ldots,\hat\e_j,\ldots,\e_n}(z,s,\ov)\right)\frac{\overline{f_{j}}}{2}dxdy\label{twointegrals}\\\nonumber+&\int_{F_M}\frac{\partial}{\partial \overline{z}}\left(\lim_{T\to\infty}D_{\e_1,\ldots,\hat\e_j,\ldots,\e_n}(z,s,\ov)\right)\frac{{f_{j}}}{2}dxdy.\end{align}
For any real-differentiable function $h:U\to\C$ where $U\subset \C$ and  any bounded domain $R\subset U$ with piecewise differentiable boundary Stokes theorem implies that \begin{equation*}2i\int_{R}\frac{\partial}{\partial \overline{z}}hdxdy=\int_{\partial R}hdz\,.\end{equation*} We apply this to the second integral in (\ref{twointegrals}).  Since $f_{j}$ is holomorphic, this integral equals
\begin{equation*}-\frac{i}{4}\int_{\partial (F_M)}\lim_{T\to\infty}D_{\e_1,\ldots,\hat\e_j,\ldots,\e_n}(z,s,\ov){f_{j}}dz.\end{equation*} The fundamental domain is the union of conjugated sides. These conjugated sides cancel in the integral. Hence this integral equals the line integral along the top of the truncated fundamental domain $F_M$. But this goes to zero by lemma \ref{stupidbound}.
We observe that when $s$ is real the first integral in (\ref{twointegrals}) is the complex conjugate of the second one. Hence this also vanishes in the limit $M\to\infty$ and we have  
 \begin{align*}\nonumber\int_{\GmodH}&\langle d\lim_{T\to\infty}D_{\e_1,\ldots,\hat\e_j,\ldots\e_n}(z,s,\ov),\alpha_{j}\rangle d\mu(z)=0.\end{align*}
We now prove that we may pull the limit outside the integral. We note that $w_j^T(z)=O(e^{-2\pi y})$ as $T\to\infty$ where the involved constant is independent of $T$. We note also that for $z\in F$ we have $\abs{-2\pi i\int_{i\infty}^{\gamma z}w_j^T}\leq\abs{\modsym{\gamma}{\a_i}}+\abs{-2\pi i\int_{i\infty}^z\a_i}$ which follows from the definition of $w_i^T$. Using this and the same approach as in the proof of lemma \ref{anotherstupidbound}, we see  that for $\Re(s)>1$ there exist $U(z,s)$ independent of $T$ such that 
\begin{align*}\abs{\modsym{dD_{\e_1,\ldots,\hat\e_j,\ldots,\e_n}(z,s,\ov)}{w_j}}\leq U(z,s)\end{align*}
and 
\begin{align*}\int_{\GmodH}U(z,s)d\mu(z)<\infty.\end{align*} Hence for any given $\varepsilon_0>0$ there exists a constant, $M$, independent of $T$ such that 
 \begin{align*}&\nonumber\abs{\int_{\GmodH}\!\!\!\left(\langle d\lim_{T\to\infty}D_{\e_1,\ldots,\hat\e_j,\ldots\e_n}(z,s,\ov),\alpha_{j}\rangle\! -\!\langle dD_{\e_1,\ldots,\hat\e_j,\ldots\e_n}(z,s,\ov),w_{j}\rangle\right) d\mu(z)}\\&
\leq\nonumber\abs{\int_{F_M}\!\!\!\left(\langle d\lim_{T\to\infty}D_{\e_1,\ldots,\hat\e_j,\ldots\e_n}(z,s,\ov),\alpha_{j}\rangle\!-\!\langle dD_{\e_1,\ldots,\hat\e_j,\ldots\e_n}(z,s,\ov),w_{j}\rangle\right) d\mu(z)}+\varepsilon_0.
\end{align*}
Hence if we choose $T>M$ and use  (\ref{getit}), (\ref{coffeeisgood}) and  Proposition \ref{approxprop} (\ref{espresso}),  we see that the integral over $F_M$ vanishes.  This finishes the proof.
\end{proof}

Using this we can now prove the following lemma: 
\begin{lemma} \label{regularfunction} The function \begin{equation}\lim_{T\to\infty}(-R(s)L_{\e_j}(\ov)D_{\e_1\ldots,\widehat{\e_j},\ldots,\e_n}(z,s,\ov))\end{equation}  is regular at $s=1$.
\end{lemma}
\begin{proof} We note that since $\a_j$ is the real part of a holomorphic differential $\delta(\a_j)=0$ and also $\delta(w_j)-\delta(\a_j)\neq 0$ only for $T \leq\Im(z)\leq T+1$ (Proposition \ref{approxprop}). We may verify that $\delta(w_j^T)=O( e^{-2\pi y})$ uniformly in $T$ and $D_{\e_1,\ldots,{\widehat{ \e_j}},\ldots,\e_n}(z,s,\ov)=O(y^{1-\sigma})$ uniformly in $T$.  Using this we  find from Lemma \ref{almostregular} that when $\Re(s)>1$
 \begin{equation}\label{hestehoved}\lim_{T\to\infty}\int_{\GmodH}L_{\e_j}(\ov)D_{\e_1,\ldots,{\widehat{ \e_j}},\ldots,\e_n}(z,s,\ov)d\mu(z)=0. \end{equation}
>From (\ref{denmarkisrainfull}) it is clear that $s=1$ is not an essential singularity. Assume that it is a pole of order $k>0.$ Hence
\begin{align}\label{contradictioncandidate}\lim_{s\to 1}(s-1)^k\lim_{T\to\infty}(-R(s)L_{\e_j}(\ov)D_{\e_1\ldots,\widehat{\e_j},\ldots,\e_n}(z,s,\ov))\neq0.\end{align}
But
\begin{align*}\lim_{s\to 1}(s&-1)^k\lim_{T\to\infty}(-R(s)L_{\e_j}(\ov)D_{\e_1\ldots,\widehat{\e_j},\ldots,\e_n}(z,s,\ov))\\&=-\lim_{s\to 1}(s-1)^k\lim_{T\to\infty}\left(\int_{\GmodH}r(z,z',s)L_{\e_j}(\ov)D_{\e_1\ldots,\widehat{\e_j},\ldots,\e_n}(z,s,\ov)\right)\\
\intertext{where $r(z,z',s)$ is the resolvent kernel}
&=-\lim_{s\to 1}\lim_{T\to\infty}\left(\int_{\GmodH}(s-1)r(z,z',s)(s-1)^{k-1}L_{\e_j}(\ov)D_{\e_1\ldots,\widehat{\e_j},\ldots,\e_n}(z,s,\ov)\right)\\
&=\vol{(\GmodH)}^{-1}\lim_{T\to\infty}\left(\int_{\GmodH}\lim_{s\to 1}(s-1)^{k-1}L_{\e_j}(\ov)D_{\e_1\ldots,\widehat{\e_j},\ldots,\e_n}(z,s,\ov)\right)\\
\intertext{since $r(z,z',s)$ has a simple pole with residue $-\vol{(\GmodH)}^{-1}$. See remark \ref{banana}}
&=\vol{(\GmodH)}^{-1}\lim_{s\to 1}(s-1)^{k-1}\lim_{T\to\infty}\int_{\GmodH}L_{\e_j}(\ov)D_{\e_1,\ldots,{\widehat{ \e_j}},\ldots,\e_n}(z,s,\ov)d\mu(z)\\
&=0
\end{align*}
by (\ref{hestehoved}). But this contradicts (\ref{contradictioncandidate}), which completes the proof.

\end{proof}
\begin{remark}\label{banana}
Using the above lemma, (\ref{denmarkisrainfull}) and the fact that the resolvent kernel  for $\Delta$, 
respectively, the Eisenstein series has expansions at 1 of the form (see e.g \cite[Theorem 2.2.6]{venkov}) 
\begin{equation}\begin{split}\label{rexpansion}r(z,z',s)&=
\frac{\vol{(\GmodH)}^{-1}}{s(1-s)}+\sum_{m=0}^\infty \widetilde{r_m}(z,z')(s-1)^m\\&
=\frac{-\vol{(\GmodH)}^{-1}}{(s-1)}+\sum_{m=0}^\infty {r_m}(z,z')(s-1)^m, \end{split}
\end{equation} 
respectively 
\begin{equation}\label{Eexpansion}E(z,s)=\frac{\vol{(\GmodH)}^{-1}}{s-1}+
\sum_{m=0}^\infty E_m(z)(s-1)^m
\end{equation} 
we may now in prin\-ciple write down the full Lau\-rent ex\-pan\-sion of the function
$\lim_{T\to\infty} D_{\e_1\ldots \e_n}(z,s,\ov)$ at $s=1$ in terms of $r_m(z,z')$,
 $E_m(z)$ and the real harmonic differentials. From this and (\ref{girlwithapearlearring}) 
we may also calculate the Laurent expansion of $E_{\e_1,\ldots,\e_n}(z,s,\ov)$ and hence of 
$E^{m,n}(z,s)$. Since general expressions are quite complicated and the combinatorics 
become quite cumbersome we restrict ourselves to some particular cases of special interest.
\end{remark}
We let $\widetilde\Sigma_{2m}$ be the elements of the symmetric group on $2m$ 
letters $1,2,\ldots,2m$ for which $\sigma(2j-1)<\sigma(2j)$ for $j=1,\ldots,m$. 
We notice that this has $(2m)!/2^m$ elements.  
\begin{lemma}\label{Dleadingterms} 
If $n$ is even $\lim_{T\to\infty}D_{\e_1,\ldots\e_n}(z,s,\ov)$ has a pole at $s=1$ of order
at most $n/2+1$. The $(s-1)^{-(n/2+1)}$ coefficient in the expansion of the function $\lim_{T\to\infty}D_{\e_1,\ldots\e_n}(z,s,\ov)$ around $s=1$ is 
\begin{equation*} 
\frac{(-8\pi^2)^{n/2}}{\vol(\GmodH)^{n/2+1}}\sum_{\sigma\in \widetilde{\Sigma}_n}
\left(\prod_{r=1}^{n/2}\int_{\GmodH}\inprod{\alpha_{\sigma(2r-1)}}{\alpha_{\sigma(2r)}}d\mu(z)\right).
\end{equation*}
If $n$ is odd,  $\lim_{T\to\infty}D_{\e_1,\ldots\e_n}(z,s,\ov)$ has a pole at $s=1$ of order at most  $(n-1)/2$.
\end{lemma} 
\begin{proof} For $n=0$ the claim is obvious, and for $n=1$ (\ref{denmarkiscold}) and  
Lemma \ref{regularfunction} give the result. Assume that the result is true 
for all $n\leq n_0$. By (\ref{denmarkisrainfull}), (\ref{Lderivedtwice}), 
Lemma \ref{regularfunction} and the fact that $\lim_{T\to\infty}(-R(s)(\modsym{w_k}{w_l})D_{\e_1,.,\hat\e_k,.,
\hat\e_l,.,\e_n}(z,s,\ov))$ can have pole order at most 1 more than $D_{\e_1,.,\hat\e_k,.,
\hat\e_l,.,\e_n}(z,s,\ov))$ at $s=1$, we obtain the result about the pole orders.
For even  $n$ we notice that by induction and using (\ref{rexpansion}) we find that the $(s-1)^{-(n/2-1)}$ coefficient is 
\begin{equation*}\begin{split}\frac{-8\pi^2}{\vol{(\GmodH)}}&\frac{(-8\pi^2)^{(n-2)/2}}{\vol(\GmodH)^{(n-2)/2+1}}\cdot\\&\sum_{\substack{k,l=1\\ k<l}}^{n}\sum_{\sigma\in \widetilde{\Sigma}_{n-2}}\left(\prod_{r=1}^{(n-2)/2}\hbox{}^\prime\int_{\GmodH}\!\!\!\!\!\inprod{\alpha_{\sigma(2r-1)}}{\alpha_{\sigma(2r)}}d\mu(z)\right)\int_{\GmodH}\!\!\!\!\!\inprod{\alpha_k}{\alpha_l}d\mu(z),\end{split}\end{equation*}
  where the prime in  the product means that we have excluded $\a_k,\a_l$ from the product and enumerated the remaining differentials accordingly. The result follows.
\end{proof}
Using this we can prove
\begin{theorem}\label{leadingterms} For all $n$  $E_{\e_1,\ldots\e_n}(z,s,\ov)$ has a pole at $s=1$ of order at most  $[n/2]+1$. If $n$ is even the  $(s-1)^{-([n/2]+1)}$ coefficient in the Laurent expansion of $E_{\e_1,\ldots\e_n}(z,s,\ov)$ is \begin{equation*} \frac{(-8\pi^2)^{n/2}}{\vol(\GmodH)^{n/2+1}}\sum_{\sigma\in \widetilde{\Sigma}_n}\left(\prod_{r=1}^{n/2}\int_{\GmodH}\inprod{\alpha_{\sigma(2r-1)}}{\alpha_{\sigma(2r)}}d\mu(z)\right).\end{equation*}
If $n$ is odd the  $(s-1)^{-([n/2]+1)}$ coefficient in the Laurent expansion of $E_{\e_1,\ldots\e_n}(z,s,\ov)$ is 
\begin{equation*}
\frac{(-8\pi^2)^{[n/2]}}{\vol(\GmodH)^{ [n/2]+1}}\sum_{k=1}^n\left(2\pi i\int_{i\infty}^z\a_k\sum_{\sigma\in\widetilde{\Sigma}_{n-1}} \prod_{r=1}^{[n/2]}\!\!\hbox{}^\prime\int_{\GmodH}\inprod{\alpha_{\sigma(2r-1)}}{\alpha_{\sigma(2r)}}d\mu(z)\right),
\end{equation*} where the prime in  the product means that we have excluded $\a_k$ from the product and enumerated the remaining differentials accordingly.
\end{theorem}
\begin{proof}This follows from (\ref{DtoE}), Lemma \ref{Dleadingterms},  and the fact that $E_{\e_1,\ldots\e_n}(z,s)$ is independent of differential within the cohomology class of the real differentials involved. \end{proof}
 We notice that \begin{equation}\label{realreal}\modsym{\Re(f(z)dz)}{\Re(f(z)dz)}=\modsym{\Im(f(z)dz)}{\Im(f(z)dz)}=y^2\abs{f(z)}^2,\end{equation} while \begin{equation}\label{realimaginary}\modsym{\Re(f(z)dz)}{\Im(f(z)dz)}=0.\end{equation} Hence many of the involved integrals may be expressed in terms of the Petersson norm defined by \begin{equation}\norm{f}:=\left(\int_{\GmodH}y^2\abs{f(z)}^2d\mu(z)\right)^{1/2}.\end{equation} We shall write $E^{\Re^l,\Im^{n-l}}(z,s):=E_{\e_1,\ldots,\e_n}(z,s,\ov)$ where $\a_i=\Re(f(z)dz)$ for $i=1,\ldots,l$ and $\a_i=\Im(f(z)dz)$ for $i=l+1,\ldots,n$. As a special case of Theorem \ref{leadingterms} we have the following theorem:
\begin{theorem}\label{moreleadingterms}
  The function $E^{\Re^{2m},\Im^{2n}}(z,s)$ has a pole of order $m+n+1$ at $s=1$, and the $(s-1)^{-(m+n+1)}$ coefficient in the Laurent expansion is \begin{equation}
\frac{{(-8\pi^2)}^{m+n}}{\vol{(\GmodH)}^{m+n+1}}\norm{f}^{2(m+n)}\frac{(2m)!(2n)!}{2^{m+n}}\binom{m+n}{n}.
\end{equation}If $n$ or $m$ are odd then the pole order of $E^{\Re^{m},\Im^{n}}(z,s)$ at $s=1$ is strictly less than $(m+n)/2+1$. 
\end{theorem}

\begin{proof}
The first part follows from Theorem \ref{leadingterms}, (\ref{realreal}) and (\ref{realimaginary}) once 
we count the number of nonzero terms in the sum indexed by $\widetilde{\Sigma}_{2m+2n}$.
 This is the set of elements  $$\widetilde{\Sigma}_{2m+2n}^{2m}=\left\{\sigma\in\widetilde{\Sigma}_{2m+2n}\left|\genfrac{}{}{0pt}{}{\sigma(2i-1),\sigma(2i)\leq 2m\textrm{ or }\sigma(2i-1),\sigma(2i)> 2m} {\textrm{ for all }i=1,\ldots, m+n} \right.\right\}.$$
This set contains $$\frac{(2m)!}{2^m}\frac{(2n)!}{2^{n}}\binom{m+n}{n}$$ elements which can be seen by noticing that each element  may be obtained uniquely by applying $\sigma_1\in \widetilde{\Sigma}_{2m}$ to $1,\ldots,2m$ and $\sigma_2\in \widetilde{\Sigma}_{2{n}}$ to $2m+1,\ldots,2m+2n$ and then shuffling $(\sigma_1(1),\sigma_1(2)),\ldots,(\sigma_1(2m-1),\sigma_1(2m))$ with $(\sigma_2(2m+1),\sigma_2(2m+2)),\ldots,(\sigma_2(2m+2n-1),\sigma_2(2m+2n))$.

If $m+n$ is odd then Theorem \ref{leadingterms} says that the pole order at $s=1$ is at most $[(m+n)/2]+1$ which is strictly less than $(m+n)/2+1$.

If $m$ and $n$ is odd then Theorem  \ref{leadingterms} says that the pole order at $s=1$ is at most $(m+n)/2+1$, but since one of the factors in the product of the $(m+n)/2+1$ term has to be zero the pole is at most of order $(m+n)/2$.
\end{proof}
We now turn to $E^{m,n}(z,s)$. We assume $f=g$.
\begin{theorem}[\cite{gos}]\label{Estar}
At $s=1$, $E^{1,0}(z,s)$ has a simple pole with residue
$$\frac{1}
{{\vol({\GmodH})}}\left(2\pi i\int_{i\infty}^zf(z)dz\right).$$
\end{theorem}
\begin{proof}
This follows directly from Theorem \ref{leadingterms} and $$E^{1,0}(z,s)=E^{\Re}(z,s)+iE^{\Im}(z,s).$$
\end{proof}

\begin{theorem}\label{thinkofanameforalabel} The Eisenstein series $E^{m,m}(z,s)$ has a pole of order $m+1$. The  $(s-1)^{-(m+1)}$ coefficient in the Laurent expansion around $s=1$ is 
\begin{equation*}\frac{(16\pi^2)^m}{\vol{(\GmodH)}^{m+1}}m!^2\norm{f}^{2m}.\end{equation*}
\end{theorem}
\begin{proof}
 Since $\modsym{\gamma}{f}=\modsym{\gamma}{\Re(f(z)dz)}+i\modsym{\gamma}{\Im(f(z)dz)}$ we have
$$\abs{\modsym{\gamma}{f}}^{2m}=(-1)^m\sum_{n=0}^m\binom{m}{n}\modsym{\gamma}{\Re(f(z)dz)}^{2n}\modsym{\gamma}{\Im(f(z)dz)}^{2(m-n)}.$$ Hence 
$$E^{m,m}(z,s)=(-1)^m\sum_{n=0}^m\binom{m}{n}E^{\Re^{2n},\Im^{2(m-n)}}(z,s).$$
 From Theorem \ref{moreleadingterms} we hence find that the leading term of $E^{m,m}(z,s)$ is \begin{equation*}\frac{(-8\pi^2)^m}{\vol{(\GmodH)}^{m+1}}\norm{f}^{2m}\sum_{n=0}^m\binom{m}{n}\frac{(2n)!(2(m-n))!}{2^m}\binom{m}{n}.\end{equation*} The sum equals $(m!)^22^m$ from which the result follows.
\end{proof}
\begin{theorem}\label{gifttogoldfeldgroup} At $s=1$, $E^{2,0}(z,s)$ has a simple pole with residue $$\frac{1}
{{\vol({\GmodH})}}\left(2\pi i\int_{i\infty}^zf(z)dz\right)^2$$ while $E^{1,1}$ has a double pole with residue
$$\frac{4\pi^2}{{\vol({\GmodH})}}\abs{\int_{i\infty}^zf(z)dz}^2+\frac{16\pi^2}{\vol}({\GmodH})
\int_{\GmodH}(E_0(z')-r_0(z,z')){y'}^2\abs{f(z')}^2d\mu(z').$$ 
The coefficient of $(s-1)^{-2}$ is
$$ \frac{16 \pi^2 \norm{f}^2}{\vol(\GmodH)^2}.$$
\end{theorem}
\begin{proof}
We start by noticing that as a special case of (\ref{DtoE}) we have  
\begin{equation*}
\begin{split}
E_{\e_1\e_2}(z, s,\ov)=&-4\pi^2\int_{i\infty}^z \a_1\int_{i\infty}^z\a_2 E(z, s)+2\pi i \int_{i\infty}
^z\a_1\lim_{T\to\infty} D_{\e_2}(z, s,\ov)
\\&+2\pi i \int_{i\infty}^z\a_2\lim_{T\to\infty} D_{\e_1}(z, s,\ov)+\lim_{T\to\infty} D_{\e_1\e_2}(z, s,\ov).
\end{split}
\end{equation*}
The first term has a simple pole at $s=1$ with residue $$\frac{-4\pi^2}{\vol{(\GmodH)}}\int_{i\infty}^z \a_1\int_{i\infty}^z\a_2,$$ while the two middle terms are regular at $s=1$ by (\ref{denmarkiscold}) and Lemma \ref{regularfunction}. The singular part of the expansion of the fourth term equals the singular part of the expansion of $$\lim_{T\to\infty}(-R(s)(L_{\e_1\e_2}E(z,s)))=8\pi^2\int_{\GmodH}r(z,z',s)\modsym{\a_1}{\a_2}E(z,s).$$ This follows from  (\ref{denmarkisrainfull}) and Lemma \ref{regularfunction}. But by using (\ref{rexpansion}) and (\ref{Eexpansion}) we find that this is
\begin{equation*}\begin{split}\frac{-8\pi^2}{\vol({\GmodH})^2}&\int_{\GmodH}\!\!\!\modsym {\a_1}{\a_2}d\mu(z)(s-1)^{-2}\\&+\frac{-8\pi^2}{\vol({\GmodH})}
\int_{\GmodH}\!\!\!\!(E_0(z')-r_0(z,z'))\modsym{\alpha_1}{\alpha_2} d\mu(z')(s-1)^{-1}.\end{split} \end{equation*} Hence we know the singular part of the expansion of $E_{\e_1,\e_2}(z,s)$ at $s=1$.
 It is easy to see that \begin{eqnarray*}E^{2,0}(z,s)=&E^{\Re^2}(z,s) +2iE^{\Re,\Im}(z,s)-E^{\Im^2}(z,s)\\E^{1,1}(z,s)=&-E^{\Re^2}(z,s) -E^{\Im^2}(z,s).\end{eqnarray*}  Using the above explicit expressions for the expansions of $E_{\e_1,\e_2}(z,s,\ov)$ now gives the result when using (\ref{realreal}) and (\ref{realimaginary}).
\end{proof}
We note that this is Theorem \ref{yetanotherstupidlabel}. We state the result for the $m+n=3$ case.
\begin{theorem} At $s=1$ $E^{3,0}(z,s)$ has a simple pole with residue $$\frac{1}
{{\vol({\GmodH})}}\left(2\pi i\int_{i\infty}^zf(z)dz\right)^3$$ while $E^{2,1}(z,s)$ has a double pole with leading term
$$\frac{32\pi^2}{\vol{(\GmodH)}^2}\left(2\pi i\int_{i\infty}^zf(z)dz\right)\norm{f}^2$$
\end{theorem}

\section{Growth on vertical lines}
By using Proposition \ref{absconvergence} we see that $E^{m,n}(z,s)=O_K(1)$ for $\Re(s)=\sigma>1$ and $z$ in a fixed compact set $K$. In this section we show that when we only require $\sigma>1/2$ then we have at most polynomial growth on the line $\Re(s)=\sigma$.

We take the opportunity to correct Theorem 1.5 in \cite{petridis}.
For simplicity assume that we have only one cusp.
We first prove : 
\begin{lemma}\label{reptilienzoo} The standard nonholomorphic Eisenstein series $E(z, s)$ 
has polynomial growth in $s$ in $\Re (s)\ge1/2$. More precisely we have for any $\varepsilon>0$ and $1/2\leq\sigma\leq 1$
\begin{equation}
E(z, \sigma+it)=O_K(\abs{t}^{1-\sigma+\epsilon})\end{equation}
for all $z\in K$, a fixed compact set in $\GmodH$.
\end{lemma}
\begin{proof}
According to \cite[p.16-17]{selberg}
the scattering function $\phi (s)$ is given by
$$\phi (s)=\frac{\sqrt{\pi}\Gamma (s-1/2)}{\Gamma (s)}ab^{1-2s}L(s),$$
where $a$, $b$ are positive constants and 
$L(s)$ is a generalized Dirichlet series with constant term $1$. In particular,
$L(s)$ tends to $1$ as $\Re (s)\rightarrow\infty$. This implies that
for $\Re (s)$ sufficiently large, say $\Re (s)\ge \sigma_0>1$, 
we have $\abs{L(s)-1}\le 1/2$. Hence 
$$\frac{E(z, s)}{L(s)}$$
is bounded for $\Re (s)\ge \sigma_0$  fixed.  By the functional equation we have
\begin{equation}\label{youngwomanwithawaterpitcher}
\begin{split}
\frac{E(z, -\sigma_0+it)}{L(-\sigma_0+it)}=&\frac{\phi (-\sigma_0+it)E(z, 1+\sigma_0-it)}{L(-\sigma_0+it)}\\=&E(z, 1+\sigma_0-it)\frac{\Gamma (-\sigma_0+it-1/2)}{\Gamma (-\sigma_0+it)}\sqrt{\pi}ab^{1+2\sigma_0-2it}.
\end{split}
\end{equation}
The asymptotics of the Gamma function (Stirling's formula) imply that the 
quotient of the Gamma factors in (\ref{youngwomanwithawaterpitcher})
is asymptotic to $\abs{t}^{-1/2}$ as $\abs{t}\rightarrow \infty$. In particular, we get
that $E(z, s)/L(s)$ is bounded on the line $\Re (s)=-\sigma_0$. We want to use Phragm\'en-Lindel\"of to conclude that it
is bounded for $-\sigma_0\le \Re (s)\le \sigma_0$, so we need to verify that  $E(z, s)/L(s)$ is of finite order in this strip.

The poles of $E(z, s)$ and $\phi (s)$ are the same.
 So $E(z, s)/L(s)$ has no poles, with the possible exception
of finite many in any vertical strip coming from
$\Gamma (s-1/2)$. These poles $\gamma_1$, $\gamma_2, \ldots \gamma_k$ can be
 easily dealt with by considering $(s-\gamma_1)(s-\gamma_2)\cdots (s-\gamma_k)
E(z, s)/L(s)$. Since $\phi (s)$ is a meromorphic function of order $\le 2$ (see \cite[Theorem 7.3]{selberg2} or \cite[Theorem 3.20]{mueller3}) and 
$\Gamma(s)$ has order $1$, we see that $L(s)$ is of finite order and by 
\cite[Th. 12.9(d) p. 164] {hejhal} or \cite[Theorem 7.3]{selberg2}  we see that $E(z, s)$ is of finite order.
We can therefore apply the Phragm\'en-Lindel\"of principle in the strip
$-\sigma_0\le \Re (s)\le \sigma_0$.

 Since $\phi (s)$ is bounded
for $\Re (s)\ge 1/2$, $\abs{\Im(s)}>1$, see \cite[Lemma 8.8]{mueller} or \cite[(8.6)]{selberg2},  we see that
$$E(z, s)=\frac{E(z, s)}{L(s)}\phi (s)\frac{\Gamma (s)}{\Gamma (s-1/2)}(\sqrt{
\pi}a)^{-1}b^{2s-1}$$
is $O(\abs{t}^{1/2})$ for $\Re (s)\ge 1/2$. 

Now we can even improve the result
by applying Phragm\'en-Lindel\"of in the strip $1/2\le \Re (s)\le 1+\delta$ for some small $\delta>0$ using the fact that $E(z,s)$ is bounded for $\Re(s)=\sigma>1$. The finite number of poles $s_0, s_1, \ldots s_k$ in this region can be dealt
by multiplying with $(s-s_0)(s-s_1)\cdots (s-s_k)$. We get as result
$$E(z, s)= O_K(\abs{t}^{1-\sigma+\epsilon})$$
for all $z\in K$, a fixed compact set in $\GmodH$.
\end{proof}
\begin{remark}\label{nopoles} We remark that  the functions $E_z(z, s)$  and $E_{\bar z}(z, s)$ have no poles in 
 $\Re (s)>1/2$, $s\notin (1/2, 1]$. This follows from \cite[Satz 10.3]{roelcke}, where this statement is proved for the Eisenstein series $E^k(z, s)$ of weight $k$.
If we set
\begin{align*}E^k(z, s)=\sum_{\gamma \in \GinfmodG}\left(\frac{c\bar z+d}{cz+d}\right)^{k/2}\Im (\gamma z)^s,\end{align*}
then $E_z(z, s)=-isE^2(z, s)/(2y)$ and $E_{\bar z}(z, s)=is E^{-2}(z, s)/(2y)$, since by term\-wise dif\-feren\-tiation we have 
\begin{eqnarray*}
E_{\bar z}(z, s)&=&\frac{is}{2}\sum_{\gamma\in\GinfmodG}\Im (\gamma z)^{s-1}\overline{(cz+d)}^{-2}\\
E_{z} (z, s)&=& \frac{-is}{2}\sum_{\gamma\in\GinfmodG}\Im (\gamma z)^{s-1}{(cz+d)}^{-2}.
\end{eqnarray*}
\end{remark}
\begin{lemma}\label{bound2}The function $D_{\e_j}(z, s,\ov)$ has polynomial growth in $s$ in $\Re (s)> 1/2$. More precisely we have for any $\varepsilon>0$ and $1/2<\sigma\leq 1$ 
\begin{equation}
  D_{\e_j}(z, \sigma+it,\ov)=O(\abs{t}^{5(1-\sigma)+\varepsilon}).\end{equation}
The constant involved depends on $\sigma$, $w_j$ and $\varepsilon$. \end{lemma}
\begin{proof}
 We have by (\ref{denmarkisrainfull}) and  (\ref{Lderivedonce})
\begin{equation}\label{girlwitharedhat}
  D_{\e_j}(z, s,\ov)=-R(s)(4\pi i \langle dE(z, s), w_j\rangle -2\pi i (\delta w_j)E(z, s)).
\end{equation}
We need to control $dE(z, s)=E_z(z, s)dz+E_{\bar z}(z, s)d\bar z$ in some sense.
Since $E_z(z, s)$
has zero Fourier coefficient   $-i(sy^{s-1}+\phi (s)(1-s)y^{-s})/2$, which could vanish only for $s=0$ or $s=1/2$, and the fact that $E_z(z, s)$ has no poles for $\Re (s)>1/2$, $s\notin (1/2, 1]$ (see Remark \ref{nopoles}), the functional equation $E_z(z, s)=\phi (s)E_z(z, 1-s)$ implies that the poles of 
 $E_z(z, s)$ are  the same  as the poles $E(z, s)$ with  the same
multiplicity with the exception of finitely many in the interval $[0, 1]$. 

We repeat the Phragm\'en-Lindel\"of argument for $E(z, s)/L(s)$  with  $E_z(z, s)/L(s)$. The only difference is that since $$E_z(z,s)=-\frac{is}{2y}E^2(z,s)$$ with  $ E^2(z,s)$ bounded on vertical lines for $\Re(s)=\sigma>1$ we start with the bound $E_z(z,s)=O(\abs{t})$ when $\Re(s)=\sigma>1$   and we get $E_z(z, s)=O(\abs{t}^{3/2})$ for $\Re(s)\geq1/2$. Applying Phragm\'en-Lindel\"of again in $1/2\leq\Re(s)\leq 1+\varepsilon$ as above we find $E_{z}(z, s)=O_K(\abs{t}^{2-\sigma+\varepsilon})$.         Similarly $E_{\bar z}(z, s)=O_K(\abs{t}^{2-\sigma +\varepsilon})$. This gives an $L^2=L^2(\GmodH,d\mu)$
bound for the function in (\ref{girlwitharedhat}) to which we apply the resolvent.
Since 
\begin{equation}\label{resolventbound}\norm{R(z)}_\infty \le \frac{1}{\hbox{dist }(z, \hbox{Spec} A)}\end{equation}
 for the resolvent of
a general self-adjoint operator $A$ on a Hilbert space and
$\hbox{dist}(s(1-s), \hbox{Spec}(\Delta) )\ge \abs{t}(2\sigma -1)$ we get for 
$1/2 <\sigma=\Re(s)\leq 1 $
$$\norm{D_{\e_j}(z, s,\ov)}_{L^2}
=O_{\sigma}\left(\frac{\abs{t}^{2-\sigma +\varepsilon}}{\abs{t}(2\sigma -1)}\right)=O_\sigma(\abs{t}^{1-\sigma +\varepsilon}).$$
We finally need to get a pointwise bound from the $L^2$ bound, for which 
we use the Sobolev embedding theorem, which in dimension $2$ implies the 
bound $\norm{u}_{\infty}\ll \norm{u}_{H^2}$, where $\norm{u}_{H^2}$ is the Sobolev $2$-norm, given by $\norm{u}_{L^2}+\norm{Pu}_{L^2}$ for any second order elliptic operator $P$.
In our case $D_{\e_j}(z,s,\ov)$ satisfies $(\Delta +s(1-s))u=-L_{\e_j}(\ov) E(z, s)$, and we get $D_{\e_j}(z,s,\ov)=O_K(\abs{t}^{3-\sigma+\varepsilon})$. We now apply Phragm\'en-Lindel\"of again in the strip $1/2+\delta\leq\Re(s)\leq 1+\delta$ for some small $\delta>0$. This gives the result. 
\end{proof}

\begin{lemma}\label{reptilienzoo2}
The function $D_{\e_1,\ldots,\e_n}(z, s,\ov)$ has polynomial growth in $t$ in $\Re (s)>1/2$. More precisely we have for any $\varepsilon>0$ and $1/2<\sigma\leq 1$ \begin{equation}
D_{\e_1,\ldots,\e_n}(z, \sigma+it,\ov)
=O (\abs{t}^{(6n-1)(1-\sigma)+\varepsilon}).\end{equation}
The involved constant depends on $\varepsilon$, $\sigma$ and $w_1,\ldots,w_n$.
\end{lemma}
\begin{proof}

This is induction in $n$. For $n=1$  we refer to Lemma \ref{bound2}. We now assume that \begin{eqnarray}
\label{fdsa}D_{\e_1,\ldots,\e_m}(z, \sigma+it,\ov)
&=&O (\abs{t}^{(6m-1)(1-\sigma)+\varepsilon})\\
\label{fdsa2}\norm{L_{\e_k}(\ov)D_{\e_1,.,\hat\e_k,.,\e_m}(z,s,\ov)}_{L^2}&=&O (\abs{t}^{(6m-1)(1-\sigma)+\varepsilon})
\end{eqnarray} whenever $m\leq n-1$.
 By (\ref{denmarkisrainfull}) we see that  we need to es\-ti\-mate the two type of terms \begin{eqnarray*}L_{\e_k\e_l}(\ov)D_{\e_1,.,\hat\e_k,.,\hat\e_l,.,\e_n}(z,s,\ov)\allowdisplaybreaks\\ L_{\e_k}(\ov)D_{\e_1,.,\hat\e_k,.,\e_n}(z,s,\ov)\end{eqnarray*} when we apply the resolvent. We can control the first  (in $L^2$) by the induction hypothesis  as we note that $L_{\e_1\e_2}(\ov)$ is a compactly supported multiplication operator (see (\ref{Lderivedtwice})). We get $$\norm{L_{\e_k\e_l}(\ov)D_{\e_1,.,\hat\e_k,.,\hat\e_l,.,\e_n}(z,s,\ov)}_{L^2}=O\left( \abs{t}^{(6(n-2)-1)(1-\sigma)+\varepsilon}\right).$$ By using that $w_i$ is compactly supported we easily deduce from (\ref{Lderivedonce}) that
\begin{align}\label{denmarkiswindy}
\nonumber\norm{L_{\e_k}(\ov)D_{\e_1,.,\hat\e_k,.,\e_n}(z,s,\ov)}_{L^2}\leq C&\left(\norm{D_{z,\e_1,.,\hat\e_k,.,\e_n}(z,s,\ov)}_{L_2(O)}\right.+\\ \norm{D_{\overline{z},\e_1,.,\hat\e_k,.,\e_n}(z,s,\ov)}&_{L_2(O)}+\left.\norm{D_{\e_1,.,\hat\e_k,.,\e_n}(z,s,\ov)}_{L_2(O)}\right), 
\end{align} where $O$ is an open set lying between the support of $w_i$ and some other compact set. We now evaluate these three terms separately. To handle the first term we note that 
\begin{equation}\label{denmarkisicy}\begin{split}\norm{D_{z,\e_1,.,\hat\e_k,.,\e_n}(z,s,\ov)}_{L_2(O)}\le& c\norm{D_{\e_1,.,\hat\e_k,.,\e_n}(z,s,\ov)} _{H^1(O)}\\\le&c\norm{D_{\e_1,.,\hat\e_k,.,\e_n}(z,s,\ov)}_{H^2(O)}\\
\leq &c'\left(\norm{D_{\e_1,.,\hat\e_k,.,\e_n}(z,s,\ov)}_{L^2(O)}\right.\\&\left.+\norm{\Delta D_{\e_1,.,\hat\e_k,.,\e_n}(z,s,\ov)}_{L^2(O)}\right).
\end{split}\end{equation} We note that by (\ref{denmarkisrainfull}) and the induction hypothesis
\begin{equation}\norm{\Delta D_{\e_1,.,\hat\e_k,.,\e_n}(z,s,\ov)}_{L^2(O)}=O(\abs{t}^{(6(n-1)-1))(1-\sigma)+\varepsilon+2}).\end{equation} Hence the left-hand side of (\ref{denmarkisicy}) is $O( \abs{t}^{(6(n-1)-1)(1-\sigma)+\varepsilon+2})$. The second term of (\ref{denmarkiswindy}) may be evaluated in the same manner, while the third term is even smaller. We thus get 
\begin{equation}\norm{L_{\e_k}(\ov)D_{\e_1,.,\hat\e_k,.,\e_n}(z,s,\ov)}_{L^2}=O(\abs{t}^{(6(n-1)-1)(1-\sigma)+\varepsilon+2})\end{equation} which certainly establishes (\ref{fdsa2}) when $m=n$. By (\ref{denmarkisrainfull}), (\ref{resolventbound}) and the above we find
\begin{align*}\norm{D_{\e_1,\ldots,\e_n}(z, s,\ov)}_{L^2}&\allowdisplaybreaks \\\le \norm{R(s)}_\infty &\left\lVert \left(\sum_{k=1}^{n} L_{\e_k}(\ov)D_{\e_1,.,\widehat{\e_k},.,\e_n}(z,s,\ov)\right.\right.\\&\left.\left.\quad+\sum_{\substack{k,l=1\\ k<l}}^{n} L_{\e_k\e_l}(\ov)D_{\e_1,.,\widehat{\e_k},.,\widehat{\e_l},.,\e_n}(z,s,\ov)\right)\right\rVert_{L^2}
\\=O(&t^{(6(n-1)-1)(1-\sigma)+\varepsilon+1}).
\end{align*}
To get a pointwise bound we also need
\begin{equation}
\norm{\Delta D_{\e_1,\ldots,\e_n}(z, s,\ov)}_{L^2}=O( \abs{t}^{(6(n-1)-1)(1-\sigma)+1+2+\varepsilon}),
\end{equation}which follows from (\ref{resolventbound}) and the above. From the Sobolev embedding theorem we get 
\begin{align*}\norm{D_{\e_1,\ldots,\e_n}(z, s,\ov)}_\infty&\leq C \norm{D_{\e_1,\ldots,\e_n}(z, s,\ov)}_{H^2}\\&=O\left(\abs{t}^{(6(n-1)-1)(1-\sigma)+3+\varepsilon}\right).\end{align*} 
We apply Phragm\'en-Lindel\"of once again in the strip $1/2+\delta\leq\Re(s)\leq1+\delta$ to finish the proof.
\end{proof}
We notice that we can get polynomial bounds on $D_{\e_1,\ldots,\e_n}(z,s)$ without using Lemma \ref{bound2}. We would then have to start the induction in  Lemma \ref{reptilienzoo2} at $n=0$ by citing Lemma \ref{reptilienzoo}. This would lead to slightly larger exponents.
Using the above lemma we conclude:
\begin{theorem}\label{growthonvertical}
The functions $E_{\e_1,\ldots,\e_n}(z, s,\ov)$ and $E^{m,n}$ have polynomial growth in $t$ in $\Re (s)\ge1/2$. More precisely we have for any $\varepsilon>0$ and $1/2<\Re(s)\leq 1$ \begin{equation}
E_{\e_1,\ldots,\e_n}(z, s,\ov)=O (\abs{t}^{(6n-1)(1-\sigma)+\varepsilon})\qquad\textrm{when  $n\geq1 $},\end{equation}
\begin{equation} E^{m,n}(z,s)
=O (\abs{t}^{(6(m+n)-1)(1-\sigma)+\varepsilon})\qquad\textrm{when  $n+m\geq1 $}.\end{equation}
The involved constants depend on $\varepsilon$, $\sigma$, $f$, $g
$ and $\a_1,\ldots,\a_n$.
\end{theorem} 
Hence we have also proved Theorem \ref{demsterbendseinebuhle}.
\section{Estimating various sums involving modular symbols}
Using the results of the previous two sections we would now like to obtain asymptotics as $T\to\infty$ for sums like  \begin{equation}\label{summatoryfunction}\sum_{\substack{\gamma\in\GinfmodG\\\Nz{\gamma}\leq T}}\omega_\gamma\end{equation}where $\omega_\gamma=1$, $\modsym{\gamma}{\alpha_1}\cdots\modsym{\gamma}{\alpha_n}$ or $\omega_\gamma=\modsym{\gamma}{f}^m \overline{\modsym{\gamma}{g}}^n$. Here $ \Nz{\gamma}=\abs{cz+d}^2$ with $c,d$ the lower row in $\gamma$ and $z\in\H$.  We let 
\begin{equation*}\tilde E(z,s)=\sum_{\gamma\in \GinfmodG}\omega_\gamma\Im(\gamma z)^s,\end{equation*} and assume that this is absolutely convergent for $\Re(s)>1$, that it  has meromorphic contination to $\Re(s)\geq h$ where $h<1$, and that as a function of $s$ it has at most polynomially growth on vertical lines. We further assume that $s=1$ is the only pole in $\Re(s)\geq h$, and that for all $\varepsilon>0$ \begin{equation}\label{omegagrowthassumption}\omega_\gamma=O(\Nz{\gamma}^\varepsilon)\textrm{ as }\Nz{\gamma}\to \infty.\end{equation} We note that  Theorem \ref{ganztreubisandasgrab}, Corollary \ref{slowgrowth} and  Theorem \ref{growthonvertical} establish these properties for the relevant Eisenstein series.

Let $\phi_U:\R\to\R$, $U\geq U_0$ ,  be a family of smooth decreasing functions with \begin{equation*}\phi_U(t)=\begin{cases}1 &\textrm{ if }t\leq1-1/U\\0 &\textrm{ if }t\geq1+1/U,  \end{cases}\end{equation*} and $\phi_U^{(j)}(t)=O(U^j)$ as $U\to\infty$. For $\Re(s)>0$ we let \begin{equation*}R_U(s)=\int_0^\infty\phi_U(t)t^{s-1}dt\end{equation*}be the Mellin transform of $\phi_U$. Then we have \begin{equation}\label{vud1}R_U(s)=\frac{1}{s}+O\left(\frac{1}{U}\right)\qquad\textrm{as }U\to\infty \end{equation} and for any $c>0$  \begin{equation}\label{vud2}R_U(s)=O\left(\frac{1}{\abs{s}}\left(\frac{U}{1+\abs{s}}\right)^c\right)\qquad\textrm{as }\abs{s}\to\infty.\end{equation} Both estimates are uniform for $\Re(s)$ bounded.   The first is a mean value estimate while the second is successive partial integration and a mean value estimate. 
Mellin inversion formula now gives
\begin{align}\sum_{\gamma\in \GinfmodG}\omega_\gamma\phi_U\left(\frac{\,\Nz{\gamma}}{T}\right)=\frac{1}{2\pi
 i}\int_{\Re(s)=2}\!\!\!\frac{\tilde E(z,s)}{y^s}R_U(s)T^sds.
\end{align} We note that by (\ref{vud2}) the integral is convergent as long as $\tilde E(z,s)$ has polynomial growth on vertical lines. We now move the line of integration to the line $\Re(s)=h$ 
 with $h<1$ by integrating along a box of some height and then letting this height go to infinity. Assuming the polynomial bounds on vertical lines 
the Phragm\'en-Lindel\"of principle implies that  there is a uniform polynomial bound $O(t^a)$ in $h\leq\Re(s)\leq 2$ (excluding a small circle around $s=1$) and using  (\ref{vud2}) we find that the contribution from the horizontal sides goes to zero. If we assume that $s=1$ is the only pole of the integrand with $\Re(s)\geq h$ then using Cauchy's residue theorem we obtain
\begin{align*}
\frac{1}{2\pi i}&\int_{\Re(s)=2}\!\!\!\frac{\tilde E(z,s)}{y^s}R_U(s)T^sds\\&=\hbox{Res}_{s=1}\left(\frac{\tilde E(z,s)}{y^s}R_U(s)T^s\right)+\frac{1}{2\pi i}\int_{\Re(s)=h}\!\!\!\frac{\tilde E(z,s)}{y^s}R_U(s)T^sds.
\end{align*}
If we choose $c=a+\varepsilon$ the last integral is $O(T^hU^{a+\varepsilon})$ uniformly for $z$ in a compact set. 

Assume that $\tilde E(z,s)$ has a pole of order $l$ with $(s-1)^{-l}$ coefficient $a_{-l}$ then, if $l>1$, we have
\begin{align}
\nonumber\hbox{Res}_{s=1}\left(\frac{\tilde E(z,s)}{y^s}R_U(s)T^s\right)&=\frac{a_{-l}}{(l-1)!y}T(\log T)^{l-1}\\\nonumber&+O(T(\log T)^{l-2}+T\log T^{l-1}/U), 
\end{align}
  since the residue divided by $T$ is a polynomial in $\log T$ of degree $l-1$ with leading coefficient  $\frac{a_{-l}}{(l-1)!y}$. This gives 
\begin{align*}  \sum_{\gamma\in \GinfmodG}\omega_\gamma\phi_U\left(\frac{\,\Nz{\gamma}}{T}\right)=& \frac{a_{-l}}{(l-1)!y}T(\log T)^{l-1}\\&+ O(T(\log T)^{l-2}+T\log T^{l-1}/U+T^hU^{a+\varepsilon}).\end{align*}
If $l=1$ then
\begin{equation*}
\hbox{Res}_{s=1}\left(\frac{\tilde E(z,s)}{y^s}R_U(s)T^s\right)=\frac{a_{-1}}{y}T+O(T/U),
\end{equation*}
and we get
\begin{equation*}
 \sum_{\gamma\in \GinfmodG}\omega_\gamma\phi_U\left(\frac{\,\Nz{\gamma}}{T}\right)= \frac{a_{-1}}{y}T+O(T/U+T^hU^{a+\varepsilon}).
\end{equation*}
 If $\tilde{E}(z,s)$ has a nonsimple pole we choose $U=\log T$ and we get
\begin{equation}\label{exp1}\sum_{\gamma\in \GinfmodG}\omega_\gamma\phi_U\left(\frac{\,\Nz{\gamma}}{T}\right)=\frac{a_{-l}}{(l-1)!y}T(\log T)^{l-1}+O(T(\log T)^{l-2}).\end{equation}
In the simple pole case we choose $U=T^{(1-h)/(a+1+\varepsilon)}$ in order to balance the error terms and we get
\begin{equation}\label{exp2}\sum_{\gamma\in \GinfmodG}\omega_\gamma\phi_U\left(\frac{\,\Nz{\gamma}}{T}\right)=\frac{a_{-1}}{y}T+O(T^{\frac{a+h+\varepsilon}{a+1+\varepsilon}}).
\end{equation}  
At this point we note that if $\omega_\gamma$ is non-negative for all $\gamma\in \GinfmodG$, then by choosing $\phi_U$ and $\tilde\phi_U$ as above and further requiring $\phi_U(t)=0$ if $t\geq 1$ and  $\tilde\phi_U(t)=1$ for $t\leq 1$, we have 
\begin{equation*}\sum_{\gamma\in\GinfmodG}\omega_\gamma\phi_U\left(\frac{\,\Nz{\gamma}}{T}\right)\leq \sum_{\substack{\gamma\in\GinfmodG\\\Nz{\gamma}\leq T}}\omega_\gamma\leq \sum_{\gamma\in\GinfmodG}\omega_\gamma\tilde\phi_U\left(\frac{\,\Nz{\gamma}}{T}\right)
\end{equation*}
from which it easily follows that the middle sum has an asymptotic expansion. As an application we use this on the usual nonholomorphic Eisenstein series and we find that 
\begin{equation}
\sum_{\substack{\gamma\in\GinfmodG\\\Nz{\gamma}\leq T}}1=\frac{T}{y\vol{(\GmodH)}}+O(T^{\frac{a+h+\varepsilon}{a+1+\varepsilon}}),
\end{equation} Now we may choose $a=1-h+\varepsilon$ (see  Lemma \ref{reptilienzoo}) and we get the following result:
\begin{lemma}\label{counting}Assume that the only pole of  $E(z,s)$ in $\Re(s)\geq h$ is $s=1$. Then
\begin{equation}\sum_{\substack{\gamma\in\GinfmodG\\\Nz{\gamma}\leq T}}1=\frac{T}{y\vol{(\GmodH)}}+O(T^{\frac{1}{2-h}+\varepsilon}).\end{equation}
\end{lemma}
Using this lemma we can now deal with the general case. 
To get a result without $\phi_U$ from (\ref{exp1}) and (\ref{exp2}) we notice that if we choose $\phi_U$ such that $\phi_U(t)=1$ for $t\leq 1$ then 
\begin{equation*}
\sum_{\gamma\in \GinfmodG}\omega_\gamma\phi_U\left(\frac{\,\Nz{\gamma}}{T}\right)=\sum_{\substack{\gamma\in\GinfmodG\\\Nz{\gamma}\leq T}}\omega_\gamma+\sum_{\substack{\gamma\in\GinfmodG\\T<\Nz{\gamma}\leq T(1+1/U)}}\omega_\gamma\phi_U\left(\frac{\,\Nz{\gamma}}{T}\right).
\end{equation*}
 Using (\ref{omegagrowthassumption}) we see that we may evaluate the last sum in the following way. For any $\varepsilon>0$ this is less than a constant times
\begin{equation*}(T(1+1/U))^\varepsilon\!\!\!\!\!\! \sum_{\substack{\gamma\in\GinfmodG\\T<\Nz{\gamma}\leq T(1+1/U)}}\!\!\!\!\!\!1 \leq 2T^\varepsilon \!\!\!\!\!\!\sum_{\substack{\gamma\in\GinfmodG\\T<\Nz{\gamma}\leq T(1+1/U)}}\!\!\!\!\!\!1.
\end{equation*} 
The sum is $O(T/U)+O(T^{\frac{1}{2-h}+\varepsilon})$ by Lemma \ref{counting}. Using this with the above choices of $U$ we get the theorem:
\begin{theorem}\label{unsmooth} If $\tilde E(z,s)$ has a pole at s=1  of order $l$ with $(s-1)^{-l}$ coefficient $a_{-l}$. If $l=1$, i.e. if the pole is simple then
\begin{equation*}
 \sum_{\substack{\gamma\in\GinfmodG\\\Nz{\gamma}\leq T}}\omega_\gamma=\frac{a_{-1}}{y}T+O(T^{\max\left(\frac{a+h}{a+1},\frac{1}{2- h}\right)+\varepsilon})
.\end{equation*}
If $l>1$ then 
\begin{equation*} \sum_{\substack{\gamma\in\GinfmodG\\\Nz{\gamma}\leq T}}\omega_\gamma=\frac{a_l}{(l-1)!y}T\log T^{l-1}+O(T\log T^{l-2})
.\end{equation*}
\end{theorem}
Using this we now get an expansion of the summatory function (\ref{summatoryfunction}) in all the cases that we studied in section \ref{polarparts}. We only state the result in a few cases.
\begin{corollary}\label{sum}Let $\alpha=\Re(f(z)dz)$ and $\beta=\Im(f(z)dz)$. Then 
\begin{align}
\label{evensum}\sum_{\substack{\gamma\in\GinfmodG\\\Nz{\gamma}\leq T}}\modsym{\gamma}{\alpha}^{2m}\modsym{\gamma}{\beta}^{2n}=&\frac{(-8\pi^2)^{m+n}\norm{f}^{2m+2n}}{y\vol(\GmodH)^{m+n+1}}&\frac{(2m)!}{m!2^{m}}\frac{(2n)!}{n!2^{n}}T\log^{m+n} T\\\nonumber&&+O(T\log^{m+n-1} T) 
 \intertext{and if $m$ or $n$ is odd then}
\label{oddsum}\sum_{\substack{\gamma\in\GinfmodG\\\Nz{\gamma}\leq T}}\modsym{\gamma}{\alpha}^{m}\modsym{\gamma}{\beta}^{n}&=O(T\log^{k} T) & 
\end{align} for some $k\in \N$ strictly less than $(m+n)/2$.
\end{corollary}
\begin{proof}
This follows from Theorem \ref{unsmooth}, Theorem \ref{growthonvertical}, Corollary \ref{slowgrowth} and Theorem \ref{moreleadingterms}, once we notice that \begin{equation}\frac{(2m)!(2n)!}{2^{m+n}(m+n)!}\binom{m+n}{n}=\frac{(2m)!}{m!2^{m}}\frac{(2n)!}{n!2^{n}}.\end{equation}
\end{proof}

\begin{corollary}\label{abssquare}We have 
\begin{equation}
\sum_{\substack{\gamma\in\GinfmodG\\\Nz{\gamma}\leq T}}\abs{\modsym{\gamma}{f}}^{2m}=\frac{(16\pi^2)^m m!}{y\vol(\GmodH)^{m+1}}\norm{f}^{2m}T\log^m T+O(T\log^{m-1} T). 
\end{equation}
\end{corollary}
\begin{proof}
This follows from Theorem \ref{unsmooth}, Theorem \ref{growthonvertical}, Corollary \ref{slowgrowth} and Theorem \ref{thinkofanameforalabel}.
\end{proof}

\begin{corollary}There exists $\delta_1>0$ such that 
\begin{equation}\sum_{\substack{\gamma\in\GinfmodG\\\Nz{\gamma}\leq T}}\modsym{\gamma}{f}=\frac{1}{y\vol(\GmodH)}\left(-2\pi i\int_{i\infty}^zf(\tau )\, d\tau\right)T+O(T^{1-\delta_1}).\end{equation}
\end{corollary}
\begin{proof}
This follows from Theorem \ref{unsmooth}, Theorem \ref{growthonvertical}, Corollary \ref{slowgrowth} and Theorem \ref{Estar}.
\end{proof}
We note that by picking $z=i$ this reproves (\ref{hochstapler}).
\begin{corollary} There exists $\delta_2>0$ such that 
\begin{equation}\label{distributionsquare}
\sum_{\substack{\gamma\in\GinfmodG\\\Nz{\gamma}\leq T}}\modsym{\gamma}{f}^2=\frac{1}{y\vol(\GmodH)}\left(-2\pi i\int_{i\infty}^zf(\tau )\, d\tau\right)^2T+O(T^{1-\delta_2}).
\end{equation}
\end{corollary}
\begin{proof}
This follows from Theorem \ref{unsmooth}, Theorem \ref{growthonvertical}, Corollary \ref{slowgrowth} and Theorem \ref{gifttogoldfeldgroup}.
\end{proof}
\begin{remark} How small we can prove $1-\delta_i$ to be in the above corollaries depends of course on how good polynomial bounds we have and how far to the left we may move the line of integration. Assuming no eigenvalues $s(1-s)\in(0,1/4)$ we can move just to the right of $s=1/2$, and using the bound of Theorem \ref{growthonvertical} we get \begin{align*}1-\delta_1&=\frac{6}{7}+\varepsilon\\
1-\delta_2&=\frac{12}{13}+\varepsilon\end{align*} for any $\varepsilon>0$. 
 \end{remark}

\section{The distribution of  modular symbols}
We now show how to obtain a distribution result for the modular symbols from the asymptotic expansions of Corollary \ref{sum}. We renormalize the modular symbols in the following way. Let \begin{eqnarray*}\label{renormalize}\widetilde{\modsym{\gamma}{f}}&=&\sqrt{\frac{\vol{(\GmodH)}}{{8\pi^2\norm{f}^2}}}\modsym{\gamma}{f}\\
\widetilde{\modsym{\gamma}{\alpha}}&=&\sqrt{\frac{\vol{(\GmodH)}}{{8\pi^2\norm{f}^2}}}\modsym{\gamma}{\alpha}\\
\widetilde{\modsym{\gamma}{\beta}}&=&\sqrt{\frac{\vol{(\GmodH)}}{{8\pi^2\norm{f}^2}}}\modsym{\gamma}{\beta}\\
 \end{eqnarray*} where  $\alpha=\Re(f(z)dz)$, $\beta=\Im(f(z)dz)$. Let furthermore \begin{equation}(\GinfmodG)^T:=\left\{\gamma\in\GinfmodG |\quad \Nz{\gamma}\leq T\right\}.\end{equation} By Lemma \ref{counting} we have \begin{equation}\label{countingexpr}\#(\GinfmodG)^T=\frac{T}{\vol({\GmodH})y}+O(T^{1-\delta}),\end{equation} for some $\delta>0$. 
Now let $X_T$ be the random variable with probability measure
\begin{equation}P(X_T\in R)=\frac{ \#\left\{\gamma\in (\GinfmodG)^T \left| \frac{\widetilde{\langle \gamma,f\rangle}}{\sqrt{\log \Nz{\gamma}}}\in R\right. \right\} }{\#(\GinfmodG)^T}.
\end{equation} for $R\subset\C$ (we set  $\widetilde{<\gamma,\alpha>}/{\sqrt{\log \Nz{\gamma}}}=\leq$ if $\norm{\gamma}_z\leq 1$. Note that there are only finitely many such elements.) We consider the moments of $X_T$
\begin{equation}M_{n,m}(X_T)=\frac{1}{\#(\GinfmodG)^T}\sum_{\gamma\in (\GinfmodG)^T}\left[\Re\left(\frac{\widetilde{\modsym{\gamma}{f}}}{\sqrt{\log \Nz{\gamma}}}\right)\right]^n\left[\Im\left(\frac{\widetilde{\modsym{\gamma}{f}}}{\sqrt{\log \Nz{\gamma}}}\right)\right]^m,\end{equation}
and note that \begin{eqnarray*}
\Re(\widetilde{\modsym{\gamma}{f}})&=&i\widetilde{\modsym{\gamma}{\beta}}\\
\Im(\widetilde{\modsym{\gamma}{f}})&=&-i\widetilde{\modsym{\gamma}{\alpha}}.\\
\end{eqnarray*}  
By partial summation we have   
\begin{align*}M_{n,m}(X_T)=\frac{i^{n+m}(-1)^m}{\#(\GinfmodG)^T}&\left(\sum_{\gamma\in(\GinfmodG)^T}\widetilde{\modsym{\gamma}{\beta}}^n\widetilde{\modsym{\gamma}{\alpha}}^m\frac{1}{\log T^{(m+n)/2}}\right.\\
&\!\!\left.+\frac{m+n}{2}\!\!\int_{c_0}^T\!\!\!\!\sum_{\gamma\in(\GinfmodG)^T}\widetilde{\modsym{\gamma}{\beta}}^n\widetilde{\modsym{\gamma}{\alpha}}^m\frac{1}{t(\log t)^{(m+n)/2+1}}dt\right),\end{align*} where $c_o=\min \{\Nz{\gamma}\, |\,\Nz{\gamma}>1\}$. If we now apply Corollary \ref{sum} and (\ref{countingexpr}) we find that as $T\to\infty$
\begin{equation}M_{n,m}(X_T)\to \begin{cases} \frac{n!}{(n/2)!2^{n/2}}\frac{m!}{(m/2)!2^{m/2}},  &\textrm{ if }m\textrm{ and }n \textrm{ are even},\\0,  &\textrm{ otherwise. }\end{cases}\end{equation}
We notice that the right-hand side is the moments of the bivariate Gaussian  distribution with correlation coefficient zero. Hence by a result due to Fr\'echet and Shohat (see \cite[11.4.C]{loeve}) we conclude the following:
\begin{theorem}\label{main}Asymptotically  $\frac{\widetilde{\langle\gamma,f\rangle}}{\sqrt{\log \Nz{\gamma}}}$ has bivariate Gaussian distribution with correlation coefficient zero. More precisely we have
\begin{equation}\frac{ \#\left\{\gamma\in (\GinfmodG)^T \left| \frac{\widetilde{\langle\gamma,f\rangle }}{\sqrt{\log \Nz{\gamma}}}\in R\right. \right\} }{\#(\GinfmodG)^T}\to \frac{1}{{2\pi}}\int_R\!\! \exp\left(-\frac{x^2+y^2}{2}\right)dxdy\end{equation} as $T\to\infty$.
\end{theorem}
As an easy corollary we get
\begin{corollary}\label{mainreal}Asymptotically  $\frac{\Re(\widetilde{\langle \gamma,f\rangle })}{\sqrt{\log \Nz{\gamma}}}$ has Gaussian distribution. More precisely we have
\begin{equation}\frac{ \#\left\{\gamma\in (\GinfmodG)^T \left| \frac{\Re(\widetilde{\langle\gamma,f\rangle )}}{\sqrt{\log \Nz{\gamma}}}\in [a,b]\right. \right\} }{\#(\GinfmodG)^T}\to \frac{1}{\sqrt{2\pi}}\int_a^b\!\! \exp\left(-\frac{x^2}{2}\right)dx\end{equation} as $T\to\infty$.
\end{corollary}
The same holds for $ \Im(\widetilde{\modsym{\gamma}{f}})$. We note that by putting $z=i$ in Corollary \ref{mainreal} and Theorem \ref{main} we obtain Theorem \ref{maintheoremforreal} and Theorem \ref{maintheorem}. 

\hbox{}

{\bf Acknowledgments}

{We would like to thank Peter Sarnak for valuable comments and
suggestions in an earlier version of our work. The second author would like to thank Roloef W. Bruggeman for valuable comments related to this work. The first arthor would like to acknowledge the financial support of the Max-Planck-Institut f\"ur Mathematik while both authors acknowledge the hospitality of the Max-Planck-Institut f\"ur Mathematik during the Activity on Analytic Number Theory, 2002. }

\end{document}